\documentclass{amsart}

\usepackage{amssymb}
\usepackage{graphicx}
\usepackage{epstopdf}
\usepackage{amscd}
\usepackage{graphics}
\newtheorem{Theorem}{Theorem}[section]
\newtheorem{Proposition}{Proposition}[section]

\def\proof{\par{\it Proof}. \ignorespaces}
\def\endproof{{\ \vbox{\hrule\hbox{%
     \vrule height1.3ex\hskip0.8ex\vrule}\hrule }}\par}

\theoremstyle{definition}
\newtheorem{Definition}[Theorem]{Definition}
\newtheorem{Example}[Theorem]{Example}

\newtheorem{Conjecture}[Theorem]{Conjecture}

\theoremstyle{remark}
\newtheorem{Remark}[Theorem]{Remark}

\numberwithin{equation}{section}

\begin{document}

\title{SINGULAR STRUCTURE OF TODA LATTICES AND 
	      COHOMOLOGY OF CERTAIN COMPACT LIE GROUPS}

\author{Luis Casian}
\address{Department of Mathematics, Ohio State University, Columbus,
OH 43210}
\email{casian@math.ohio-state.edu}

\author{Yuji Kodama$^*$}
\thanks{$^*$Partially
supported by NSF grant DMS0404931}

\address{Department of Mathematics, Ohio State University, Columbus,
OH 43210}
\email{kodama@math.ohio-state.edu}



\thispagestyle{empty}
\pagenumbering{roman}

\pagenumbering{arabic}
\setcounter{page}{1}

\begin{abstract}
 We study the singularities (blow-ups) of the Toda lattice associated with a real split semisimple Lie algebra $\mathfrak g$. It turns out that 
the total number of blow-up points along trajectories of the Toda lattice is given by the number of points of a Chevalley group $K({\mathbb F}_q)$ related to the maximal compact subgroup $K$
of the group $\check G$ with $\check{\mathfrak g}={\rm Lie}(\check G)$ over the finite field ${\mathbb F}_q$. Here $\check{\mathfrak g}$ is the Langlands dual of ${\mathfrak g}$.
The blow-ups of the Toda lattice are given by the zero set of the $\tau$-functions. For example, the blow-ups of the Toda lattice of A-type are determined by the zeros of the Schur polynomials associated with rectangular Young diagrams. Those Schur polynomials are the $\tau$-functions for the nilpotent Toda lattices. Then we conjecture that the number of
blow-ups is also given by the number of real roots of those Schur polynomials for a specific variable.
We also discuss the case of periodic Toda lattice in connection with the real cohomology of
the flag manifold associated to an affine Kac-Moody algebra.
\end{abstract}

\maketitle

\section{Introduction: Toda lattices and the blow-ups}
Let us first give some notations and definitions of the 
 real split semisimple Lie algebra ${\mathfrak g}$ of rank $l$: We fix a split Cartan subalgebra ${\mathfrak h}$ with root system $\Delta=\Delta({\mathfrak g},{\mathfrak h})=\Delta^+\cup\Delta^-$, real root vectors $e_{\alpha_i}$ associated with simple roots $\Pi=\{\alpha_i:i=1,\ldots,l\}$. We also denote $\{h_{\alpha_i},e_{\pm\alpha_i}\}$ the Cartan-Chevalley basis of the algebra ${\mathfrak g}$ which satisfies the relations,
\[
    [h_{\alpha_i} , h_{\alpha_j}] = 0, \quad
    [h_{\alpha_i}, e_{\pm \alpha_j}] = \pm C_{j,i}e_{\pm \alpha_j} \ , \quad
    [e_{\alpha_i} , e_{-\alpha_j}] = \delta_{i,j}h_{\alpha_j},
\]
where $(C_{i,j})$ is the $l\times l$ Cartan matrix of the Lie algebra
$\mathfrak g$ and $C_{i,j}=\alpha_i(h_{\alpha_j})$ (as used in \cite{humphreys:72}). 
For example, the Cartan matrices for $B_2,C_2$ and $G_2$ are given by,
\[
B_2:~
\begin{pmatrix} 2&-2\\-1&2\end{pmatrix}
\quad C_2:~ \begin{pmatrix}2&-1\\-2&2\end{pmatrix}
\quad G_2:~
\begin{pmatrix}2&-1\\-3&2\end{pmatrix}
\]
The Lie algebra ${\mathfrak g}$ admits the decomposition,
\[
{\mathfrak g}={\mathfrak n}^-\oplus{\mathfrak h}\oplus{\mathfrak n}^+={\mathfrak n}^-\oplus {\mathfrak b}^+={\mathfrak b}^-\oplus{\mathfrak n}^+,
\]
where ${\mathfrak n}^{\pm}$ are nilpotent subalgebras defined by ${\mathfrak n}^{\pm}=\oplus_{\alpha\in\Delta^{\pm}}{\mathbb R}e_{\alpha}$ with root vectors $e_{\alpha}$, and ${\mathfrak b}^{\pm}={\mathfrak n}^{\pm}\oplus{\mathfrak h}$ are Borel subalgebras of ${\mathfrak g}$.
We denote by  $\check {\mathfrak g}$ the  real split Lie algebra with Cartan matrix given by  the transpose of the Cartan matrix of $\mathfrak g$ ($\check{\mathfrak g}$
is called the Langlands dual of ${\mathfrak g}$). Note that ${\mathfrak g}=\check {\mathfrak g}$ if ${\mathfrak g}$ is simple and not of type $B$ or $C$.  

We also fix a split Cartan subgroup $H$ with ${\rm Lie}(H)={\mathfrak h}$ and a Borel subgroup $B$
with ${\rm Lie}(B)={\mathfrak b}^+$ with $B=HN$ where $N$ is a Lie group having the Lie algebra ${\mathfrak n}^+$. We also denote Lie groups $B^-$ and $N^-$ with ${\rm Lie}(B^-)={\mathfrak b}^-$
and ${\rm Lie}(N^-)={\mathfrak n}^-$.
Integral weights on ${\mathfrak h}$ can be exponentiated to $H$. For example if
$\alpha$ is a root then there is a corresponding character $\chi_\alpha$ defined on $H$.

\medskip

Most of the results presented in this paper can be found in our recent paper \cite{casian:05},
and the main purpose of this paper is to give a brief summary of those, putting emphasis on the singular structure of the Toda lattice.
In addition, we will also discuss an extension to the case of periodic Toda lattice whose underlying algebra is given by an affine Kac-Moody algebra.

\subsection{Toda lattices: nonperiodic case}
The nonperiodic Toda lattice equation related to
real split semisimple Lie algebra $\mathfrak g$ of rank $l$ is defined by
the Lax equation,
(\cite{bogoyavlensky:76,kostant:79}),
\begin{equation}
\label{lax}
\displaystyle{\frac{dL}{dt}=[L,A]}
\end{equation}
where $L$ is a Jacobi element of ${\mathfrak g}$ and $A$ is the ${\mathfrak n}^-$-projection
of $L$, denoted by $\Pi_{\mathfrak n^-}L$,
\begin{equation}\left\{
\begin{array} {ll}
& \displaystyle{L(t)=\sum_{i=1}^l b_i(t)h_{\alpha_i}+\sum_{i=1}^l \left(
a_i(t)e_{-\alpha_i}+e_{\alpha_i}\right)} \\
& \displaystyle{A(t)=\Pi_{\mathfrak n^-}L=\sum_{i=1}^l a_i(t)e_{-\alpha_i}}
\end{array}
\right.
\label{LA}
\end{equation}
The Lax equation (\ref{lax}) then gives
\begin{equation}
\left\{
\begin{array} {ll}
& \displaystyle{\frac{d b_i}{dt}=a_i\,,} \\
& \displaystyle{\frac{d a_i}{dt}=  -\left(\sum_{j=1}^l C_{i,j} b_j\right)a_i\,.}
\end{array}
\right.
\label{toda-lax}
\end{equation}

The integrability of the system can be shown by the existence of the Chevalley
invariants, $\{I_k(L): k=1,\cdots,l\}$, which are given by the homogeneous
polynomial
of $\{(a_i, b_i) : i=1,\cdots,l\}$. (Recall that those correspond to the basic invariants in ${\mathbb C}[{\mathfrak h}]$ of the Weyl group $W$, i.e.
 ${\mathbb C}[{\mathfrak g}]^G\cong {\mathbb C}[{\mathfrak h}]^W$ with Ad-action of $G$.)
 The invariant polynomials also define the commutative equations
of the Toda equation (\ref{lax}),
\begin{equation}
\label{higherflows}
\frac{\partial L}{\partial t_k}=[L,\Pi_{\mathfrak n^-}\nabla I_k(L)]\,
\quad {\rm for}\quad k=1,\cdots,l\,,
\end{equation}
where $\nabla$ is the gradient with respect to the Killing form, i.e.
for any $x\in {\mathfrak g}$, $dI_k(L)(x)=K(\nabla I_k(L),x)$. Here $\{t_k: k=1,\ldots,l\}$ represent the flow
parameters, and we will also denote $t_k$ by $t_{m_k}$ with the exponent $m_k$ of the basic invariant $I_k$ (recall $m_k=d_k-1$ where $d_k$ is the degree of $I_k$, see e.g. \cite{carter:89}).
For example, in the case of ${\mathfrak g}={\mathfrak{sl}}(l+1;{\mathbb R})$,
the invariants $I_k(L)$ and the gradients $\nabla I_k(L)$ are given by
\[I_k(L)=\frac{1}{k+1}{\rm tr}(L^{k+1})\,\quad{\rm and}\quad
\nabla I_k(L)=L^k.\]
In this case, the degree of $I_k$ is $k+1$ and the exponent is $m_k=k$.
The set of commutative equations (\ref{higherflows}) is called the Toda lattice hierarchy.
Note that Eq. (\ref{toda-lax}) is the first member of the hierarchy, i.e. $t=t_1$.
Then the {\it real} isospectral
manifold is defined by
\[
Z(\gamma)_{\mathbb R}=\left\{(a_1,\cdots,a_l,b_1\cdots,b_l)\in{\mathbb
R}^{2l}~:~
I_k(L)=\gamma_k\in {\mathbb R},~ k=1,\cdots,l\right\}.
\]
The manifold $Z(\gamma)_{\mathbb R}$ can be compactified by adding the set of
points corresponding to the singularities ({\it blow-ups}) of the solution.
Then the compact manifold ${\tilde Z}(\gamma)_{\mathbb R}$ is
described by a union of convex
polytopes $\Gamma_{\epsilon}$ with $\epsilon=(\epsilon_1,\cdots,\epsilon_l),~
\epsilon_i={\rm sgn}(a_i)$ \cite{casian:02},
\[
{\tilde Z}(\gamma)_{\mathbb R}=\bigcup_{\epsilon\in\{\pm\}^l}\Gamma_{\epsilon}.
\] 
Each polytope
$\Gamma_{\epsilon}$ is expressed as the closure of the orbit of a
Cartan subgroup. Thus in an ad-diagonalizable
case with distinct eigenvalues, the
compact manifold ${\tilde Z}(\gamma)_{\mathbb R}$ is a toric variety, and
the vertices of each polytope are labeled by the elements of the Weyl group.

\begin{Example}: For the case ${\mathfrak g}={\mathfrak{sl}}(l+1;{\mathbb R})$, i.e. $A_l$-type, we have
\[
L=\begin{pmatrix} 
b_1   &   1    &     0     &     \cdots      &     0    \\
a_1   &  b_2-b_1 &   1     &     \cdots      &     0   \\
0        &  a_2 &  b_3-b_2  &    \cdots      &     0    \\
\vdots &  \vdots  & \vdots  &  \ddots  &  \vdots \\
0     &      0       &   0       &   \cdots    &   -b_{l}
\end{pmatrix}
\]
Notice that a fixed point $(a_1=\cdots=a_l=0)$ is a triangular matrix with
the eigenvalues on the diagonal, and the total number of fixed points is given by $(l+1)! =|S_{l+1}|$ (assuming all eigenvalues are distinct).
Each fixed point is then labeled by a unique element of the Weyl group $S_{l+1}$, and
it is identified as a vertex of the isospectral polytope $\Gamma_{\epsilon}$.
Then the isospectral polytope is just a permutohedron associated with $S_{l+1}$.
For example, $A_2$-Toda lattice, we have
a hexagon as the isospectral polytope whose vertices are the fixed points
with the Lax matrices given by
\[
L_w:=\begin{pmatrix}
\lambda_{w^{-1}(3)}  &   1   &  0   \\
0   & \lambda_{w^{-1}(2)} & 1    \\
0  &  0  &  \lambda_{w^{-1}(1)}
\end{pmatrix}\quad {\rm with}\quad\lambda_1>\lambda_2>\lambda_3
\]
where $w^{-1}(i)$ represents a permutation of $(3,2,1)$ associated to $w\in S_3$
(see Figure \ref{A2:fig}). In Figure \ref{A3:fig}, we show the polytope $\Gamma_{\epsilon}$
for $A_3$ Toda lattice. The vertices of the polytope are marked by the elements of the
Weyl group $S_4$. The faces and edges on the boundary of the polytope
correspond to the subsystems given by $a_k=0$ for some $k$.
\begin{figure}[t]
\includegraphics[width=5cm]{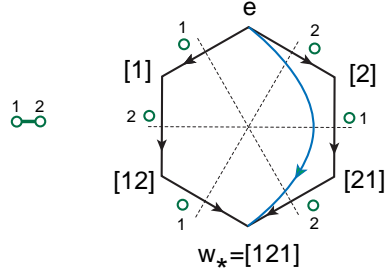}
\caption{Isospectral polytope $\Gamma_{\epsilon}$ of $A_2$-Toda lattice. The Dynkin diagram of $A_2$ is shown on the left, and each edge corresponds to $A_1$-Toda lattice
whose Dynkin diagram is just one circle. For example, the edge from $e$ to $[1]=s_1$ corresponds to the system with $a_2=0$.  The element $e\in S_3$ corresponds to the
Lax matrix $L$ whose diagonal part is given by $(\lambda_3,\lambda_2,\lambda_1)$.
Then the element $[1]$ describes the Lax matrix with diag$(\lambda_2,\lambda_3,\lambda_1)$,
i.e. we have the $s_1$-action on the diagonal of $L$. Each element $w$ then
corresponds to the Lax matrix with $w^{-1}\cdot (3,2,1)=(w^{-1}(3),w^{-1}(2),w^{-1}(1))$.}
\label{A2:fig}
\end{figure}
In particular, each edge of the isospectral polytope corresponds to a $A_1$-Toda lattice,
which is given by the Lax pair,
\[
L=\begin{pmatrix} \tilde b_k & 1 \\ a_k& \tilde b_{k+1}\end{pmatrix}\,,\quad\quad
A=\begin{pmatrix} 0 & 0 \\ a_k & 0 \end{pmatrix}\,.
\]
If $a_k(0)>0$, the flow is complete, and $a_1(t_1)\to 0$ as $t_1\to\pm\infty$.
The isospectral polytope is a connected line segment, denoted by $\Gamma_+$, with
the end-points corresponding to the matrices with $\lambda>\mu$,
\[
L(t_1=-\infty)=\begin{pmatrix} \mu & 1 \\ 0 & \lambda \end{pmatrix}\quad
L(t_1=\infty)=\begin{pmatrix} \lambda & 1 \\ 0 &\mu \end{pmatrix}
\]
 If $a_k(0)<0$, the flow has a singularity (blow-up) in finite time.
The isospectral polytope consists of two line segments, denoted by $\Gamma_-$
(see Figure \ref{A1:fig}).

\begin{figure}[t]
\includegraphics[width=7cm]{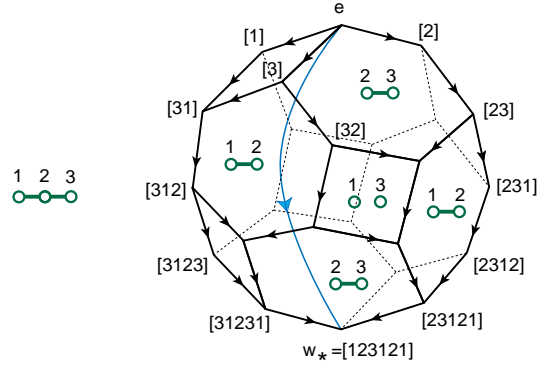}
\caption{Isospectral polytope $\Gamma_{\epsilon}$ for $A_3$-Toda lattice. The boundary of the polytope consists of the subsystems of the Toda lattice, in which eight hexagons
corresponds to $A_2$-Toda lattices and six squares corresponds to $A_1\times A_1$-Toda lattices.
Those subsystems are marked by the corresponding sub-Dynkin diagrams.}
\label{A3:fig}
\end{figure}

\end{Example}

\subsection{The $\tau$-functions and Painlev\'e divisors}
The analytical structure of the blow-ups can be obtained by the
$\tau$-functions,
which are defined by
\begin{equation}
\label{ab}
b_k=\displaystyle{\frac{d}{dt_1}\ln \tau_k, \quad \quad a_k=a_k^0\prod_{j=1}^l
(\tau_j)^{-C_{k,j}}},
\end{equation}
where $a_k^0$ are some constants.
  The tau-functions are given by \cite{flaschka:91},
\begin{equation}
\label{tau}
\tau_j(t_1,\ldots,t_l)=\langle g(t_1,\ldots,t_l)\cdot v^{\omega_j},v^{\omega_j}\rangle,\quad
g=\exp\left(\sum_{k=1}^l t_k\nabla I_k(L^0)\right)\,.
\end{equation}
Here $v^{\omega_j}$ is the highest weight vector in the fundamental representation
of $G$, and $\langle\cdot,\cdot\rangle$ is a pairing on the representation space,
and $L^0$ is an initial data of $L(t_1,\ldots,t_l)$.
The blow-up points (i.e. the singular points of $(a_j,b_j)$) are given by the zeros of the $\tau$-functions, $\tau_j(t_1,\ldots,t_l)=0$
for some $j\in\{1,\ldots,l\}$. We then define the Painlev\'e divisor ${\mathcal D}_J$ for
a subset $J\subset\{1,\ldots,l\}$ as (\cite{casian:02})
\[
{\mathcal D}_J:=\bigcap_{j\in J}\{\tau_j(t_1,\ldots,t_l)=0\},\quad {\rm codim}_{\mathbb R}{\mathcal D}_J=
|J|.
\]
The ${\mathcal D}_J$ can be also described by the intersection with the Bruhat cell
$N^-w_JB/B$ with the longest element $w_J$ of the Weyl subgroup $W_J=\langle
s_j : j\in J\rangle$ (Theorem 3.3 in \cite{flaschka:91}).
In particular, the divisor ${\mathcal D}_{\{1,\ldots,l\}}$ is a unique point, denoted as $p_o$, in the variety
$\tilde Z(\gamma)_{\mathbb R}$, and it is contained in the $\Gamma_{\epsilon}$-polytope with $\epsilon=(-\cdots -)$. 
Then the geometry of the divisor ${\mathcal D}_0=\cup_{j=1}^l \{\tau_j=0\}$, the union of the Painlev\'e
divisors ${\mathcal D}_{\{j\}}$, near the point $p_o$ can be expressed as the product of
$\tau$-functions,
\begin{equation}
\label{divF}
F(t_1,\ldots,t_l):=\prod_{j=1}^l\tau_j(t_1,\ldots,t_l)=F_d(t_1,\ldots,t_l)+F_{d+1}(t_1,\ldots,t_l)+\cdots\,,
\end{equation}
where each $F_k(t_1,\ldots,t_l)$ is a homogeneous polynomial of degree $k$.
The algebraic variety $V:=\{F_d=0\}$ defines the tangent cone at the point $p_o$, and the degree
$d$ is the {\it multiplicity} of the singularity of $V$ at $p_o$.
The number $d$ has several surprising connections with other numbers, such as the number of ${\mathbb F}_q$ points on the maximal compact subgroup of the underlying group of
the Toda lattice and the number of real roots of certain symmetric functions (e.g. Schur polynomials).
Here ${\mathbb F}_q$ is a finite field with $q$ elements, with $q$  a power of a prime.
One of the main purpose of this paper is to explain those connections (the details can be found in
our recent paper \cite{casian:05}).

\subsection{Action of the Weyl group on the signs of the Toda lattice.}
Here we give an algebraic description of the blow-ups, so that one can compute
the number of blow-ups in the Toda flow. 
The following action of the Weyl group $W$ describes how the signs of the functions
$a_j$ for $j=1,\cdots, l$ change when $a_i$ blows up.

\begin{Definition} \label{act} (Proposition 3.16 in \cite{casian:02}) For any set of signs
$\epsilon =(\epsilon_1,\cdots , \epsilon_l)\in \{\pm\}^l$,
a simple reflection $s_i:=s_{\alpha_i}\in W$ acts on the sign $\epsilon_j$ by
\begin{equation}
\nonumber
s_i~:~\epsilon_j \longmapsto \epsilon_j\epsilon_i^{-C_{j,i}}.
\label{Waction}
\end{equation}
The sign change is defined on the group character $\chi_{\alpha_i}$ with
$\epsilon_i={\rm sign}(\chi_{\alpha_i})$ (recall $s_i\cdot\alpha_j=\alpha_j-C_{j,i}\alpha_i)$. 
We also identify the sign $\epsilon_i$
as that of $a_i$, since the functions $a_k$ are given by (\ref{ab}) and the $\tau$-functions are
given by the fundamental weights $\omega_j$ in (\ref{tau}), which relate to the equation
$\chi_{\alpha_j}=\prod_{k=1}^l(\chi_{\omega_k})^{C_{j,k}}$ (recall $\alpha_j=\sum_{k=1}^l C_{j,k}\omega_k$).

We now define the relation $\Rightarrow$ between the vertices of the polytope
$\Gamma_{\epsilon}$ as follows: if  $\epsilon_i=+$ then we  write: 
$\epsilon \overset{s_i}{\Rightarrow}\epsilon^\prime$ where  $\epsilon^\prime=s_i\epsilon$. We
also write $w\Rightarrow ws_i$ (see Example \ref{A1:eg} below).
\end{Definition}

Under the action of the Weyl group, not every simple reflection $s_i$ changes the
sign $\epsilon$. The following is an alternative way to measure the size of $w$
which only takes into account simple reflections that change the sign $\epsilon$,
that is, a trajectory of a Toda lattice having a blow-up point.
These numbers will later reappear in the context of the computation of certain
Frobenius eigenvalues.

Now the following definition gives the number of blow-ups in the Toda orbit from
the top vertex $e$ to the vertex labeled by $w\in W$:
\begin{Definition} \label{eta} Choose a reduced expression $w=s_{j_1}\cdots s_{j_r}$.  Consider the sequence of signs as the orbit given by $w$-action:
\[ \epsilon~ \to ~s_{j_1}\epsilon~ \to ~s_{j_2} s_{j_1} \epsilon ~\to~ \cdots~\to~w^{-1}\epsilon\,.
\]
We then define the function $\eta(w,\epsilon)$ as the number of $\to$ which are not of the form $\overset{s_i}{\Rightarrow}$ as in Definition \ref{act}. The number $\eta(w_*,\epsilon)$
for the longest element $w_*$ gives the total number of blow-ups along the Toda flow
in $\Gamma_{\epsilon}$-polytope. Whenever $\epsilon=(-\cdots -)$ we 
will just denote $\eta(w,\epsilon)=\eta(w)$.
\end{Definition}
Note that each reduced expression of $w$ corresponds to a path following  Toda lattice
trajectories along 1-dimensional subsystems leading to $w$.  Each 1-dimensional
subsystem is equivalent to $A_1$-Toda lattice.
In Corollary 5.2 of \cite{casian:05} it is shown that
$\eta(w,\epsilon)$ is independent of the reduced expression.  Hence the number of blow-up
points along trajectories in one-dimensional subsystems in the boundary of
the $\Gamma_{\epsilon}$-polytope is independent of the trajectory (parametrized
by a reduced expression).  
\begin{figure}[t]
\includegraphics[width=3.7in]{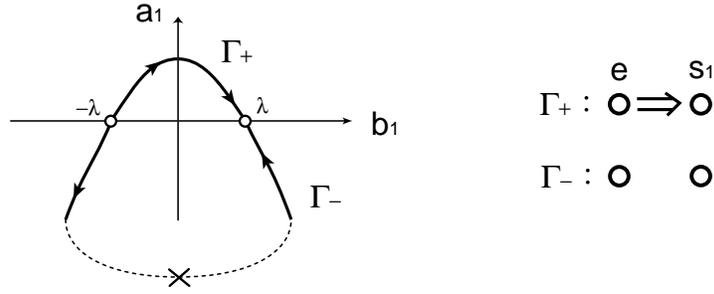}
\caption{The isospectral manifold for the $A_1$-Toda lattice.
The left figure shows the invariant curve $I_1=a_1+b_1^2=\lambda^2$,
and the mark $\times$ indicates the blow-up point, i.e. the divisor ${\mathcal D}_{\{1\}}$. The right figure shows the
graphs of $\Gamma_{\pm}$-polytopes with the Weyl action on the signs.
This shows that $\Gamma_+$ is connected, and $\Gamma_-$ has two connected components
with ${\mathcal D}_{\{1\}}$. }
\label{A1:fig}
\end{figure}

\begin{Example}\label{A1:eg}
We consider the ${\mathfrak{sl}}(2;{\mathbb R})$-Toda lattice which is the simplest case, but
provides the basic structure of the general case.
The Lax pair $(L,A)$ is given by
\[
L=\begin{pmatrix}
b_1 & 1 \\ a_1  &  -b_1\end{pmatrix},
\quad A=\begin{pmatrix} 0  &  0 \\ a_1 &  0  \end{pmatrix}.
\]
The Chevalley invariant is given by $I_1=\frac{1}{2}{\rm Tr}(L^2)=a_1+b_1^2$.
Then the isospectral manifold $Z(\gamma)_{\mathbb R}$ for a real split case is given by
the curve $I(a_1,b_1)=\gamma_1>0$ (see Figure \ref{A1:fig} where $\gamma_1=\lambda^2$). The compactified manifold
${\tilde{Z}}(\gamma)_{\mathbb R}$ consists of two polytopes (line segments) $\Gamma_+$
and $\Gamma_-$, 
\[
{\tilde Z}(\gamma)_{\mathbb R}=\Gamma_+\cup \Gamma_- \cong S^1\,.
\]
Here the $\Gamma_-$ is compactified by adding the blow-up point marked by
$\times$ in Figure, i.e. $\tau_1=0$.
The end-points (vertices) of each segment are marked by the Weyl elements, $e$ and $s_1$.
Figure \ref{A1:fig} also shows the graphs associated with the polytopes $\Gamma_{\pm}$
where the connection with the arrow in $\Gamma_+$ indicate no blow-up in the flow between
the vertices (see Definition \ref{graph} for more details).

The solution $(a_1,b_1)$ can be expressed by the $\tau$-function:
\begin{itemize}
\item[a)] For the case $a_1>0$ (i.e. $\Gamma_+$-polytope), Eq. (\ref{tau}) gives
\[
\tau_1(t_1)={\rm cosh}(\lambda t_1),
\]
which leads to the solution
\[ a_1(t_1)=\lambda^2 {\rm sech}^2(\lambda t_1),\quad b(t_1)=\lambda {\rm tanh}(\lambda t_1).
\]
Since there is no blow-up in this case, we have $\eta(e,+)=\eta(s_1,+)=0$. This implies
that there is an edge between the vertices in $\Gamma_+$ (see Figure \ref{A1:fig}).
\item[b)] For the case $a_1<0$ (i.e. $\Gamma_-$-polytope), we have
\[
\tau_1(t_1)=\frac{1}{\lambda}{\rm sinh}(\lambda t_1),
\]
which gives
\[
a_1(t_1)=-\lambda^2{\rm csh}^2(\lambda t_1),\quad b_1(t_1)=\lambda {\rm coth}(\lambda t_1).
\]
Thus the solution $(a_1(t_1), b_1(t_1))$  blows up at $t_1=0$. We have $\eta(e)=0$ and $\eta(s_1)=1$,
which gives no connecting arrow between the vertices in $\Gamma_-$ as in Figure \ref{A1:fig}.
\end{itemize}
Note that in both cases the solution approaches  the fixed points $(a_1=0,b_1=\pm\lambda)$
as $t\to\pm\infty$, which are the vertices of the polytope.
\end{Example}

\begin{figure}[t]
\includegraphics[width=4.7in]{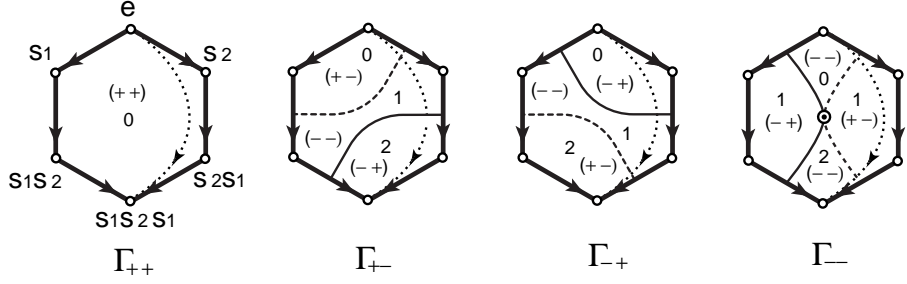}
\caption{The four hexagons associated to the signs $\epsilon=(\epsilon_1,\epsilon_2)$. 
The Painlev\'e divisors ${\mathcal D}_{\{1\}}$ and ${\mathcal D}_{\{2\}}$
 are indicated by the solid curves for the blow-ups of $a_1$ (i.e. $\tau_1=0$)
and the dashed curves for the blow-up of  $a_2$ (i.e. $\tau_2=0$).
The double circle in $\Gamma_-$ indicate the point $p_o$ corresponding to
$\cap_{j=1}^2{\mathcal D}_{\{j\}}$.
The boundaries of the hexagons describe the subsystems given by $a_i=0$ for $i=1,2$.
 The compactification can be done uniquely by gluing the boundaries
according to the sign changes. Each hexagon is divided by the 
Painlev\'e divisors into connected components. A Toda flow in $t_1$-variable
is shown as the dotted curve starting from
the vertex marked by the identity
element $e$, and ending to the vertex by the longest element $w_*=s_1s_2s_1$.}
\label{hexagon1:fig}
\end{figure}

\begin{figure}[t]
\includegraphics[width=5.2in]{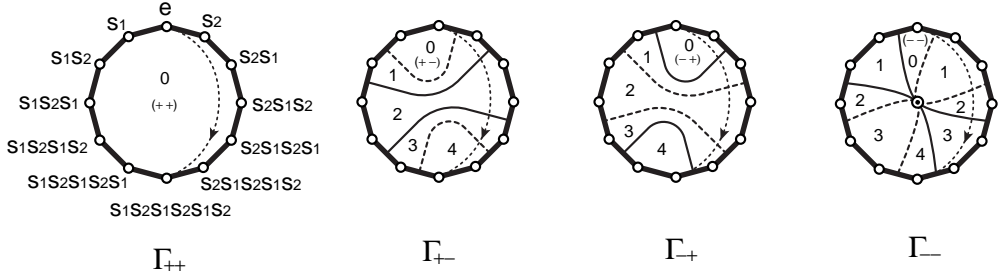}
\caption{The 12-gons corresponding to signed Toda lattice for $G_2$.
The Painlev\'e divisors ${\mathcal D}_{\{1\}}$ and ${\mathcal D}_{\{2\}}$ corresponding to the blow ups of $a_1$ and $a_2$ are 
indicated by the solid and dashed curves
inside the 12-gons. A Toda flow in $t_1$-variable
is shown as the dotted curve starting from
the vertex marked by the identity
element $e$, and ending to the vertex by the longest element $w_*$.}
\label{G2:fig}
\end{figure}

\begin{Example}
The cases of $A_2$ and $G_2$ are illustrated in Figures \ref{hexagon1:fig} and \ref{G2:fig}.
In these figures the four hexagons and 12-gons are shown as the $\epsilon$-polytope
$\Gamma_{\epsilon}$ with the signs $\epsilon=(\epsilon_1\epsilon_2)$.
These polytopes glue together to form a compact isospectral manifold $\tilde Z(\gamma)_{\mathbb R}$.
Trajectories of the Toda lattice starts in the vertex associated to $e$ and move towards the vertex corresponding to the longest element in the Weyl group. 

In the case of $A_2$, the $W$-action on the signs $\epsilon=(\epsilon_1\epsilon_2)$ gives
  $s_1(--)=(-+)$,
$s_2(-+)=(-+)$ and $s_1(-+)=(--)$. From those we obtain $\eta (e)=0$, $\eta(s_1)=\eta(s_1s_2)=1$, $\eta(s_2)=\eta(s_2s_1)=1$ and $\eta(s_1s_2s_1)=2$. Those give the numbers of blow-ups
in the Toda flow (see Figure \ref{hexagon1:fig}).

In the case of $G_2$, we obtain $\eta(e)=0$, $\eta(s_1)=\eta(s_2)=\eta(s_1s_2)=\eta(s_2s_1)=1$,
$\eta(s_1s_2s_1)=\eta(s_2s_1s_2)=2$, $\eta(s_1s_2s_1s_2)=\eta(s_2s_1s_2s_1)=\eta(s_1s_2s_1s_2s_1)=\eta(s_2s_1s_2s_1s_2)=3$
and $\eta(w_*)=4$. The total number of blow-ups is then 4 (see Figure \ref{G2:fig}).
\end{Example}

\begin{figure}[t]
\includegraphics[width=9cm]{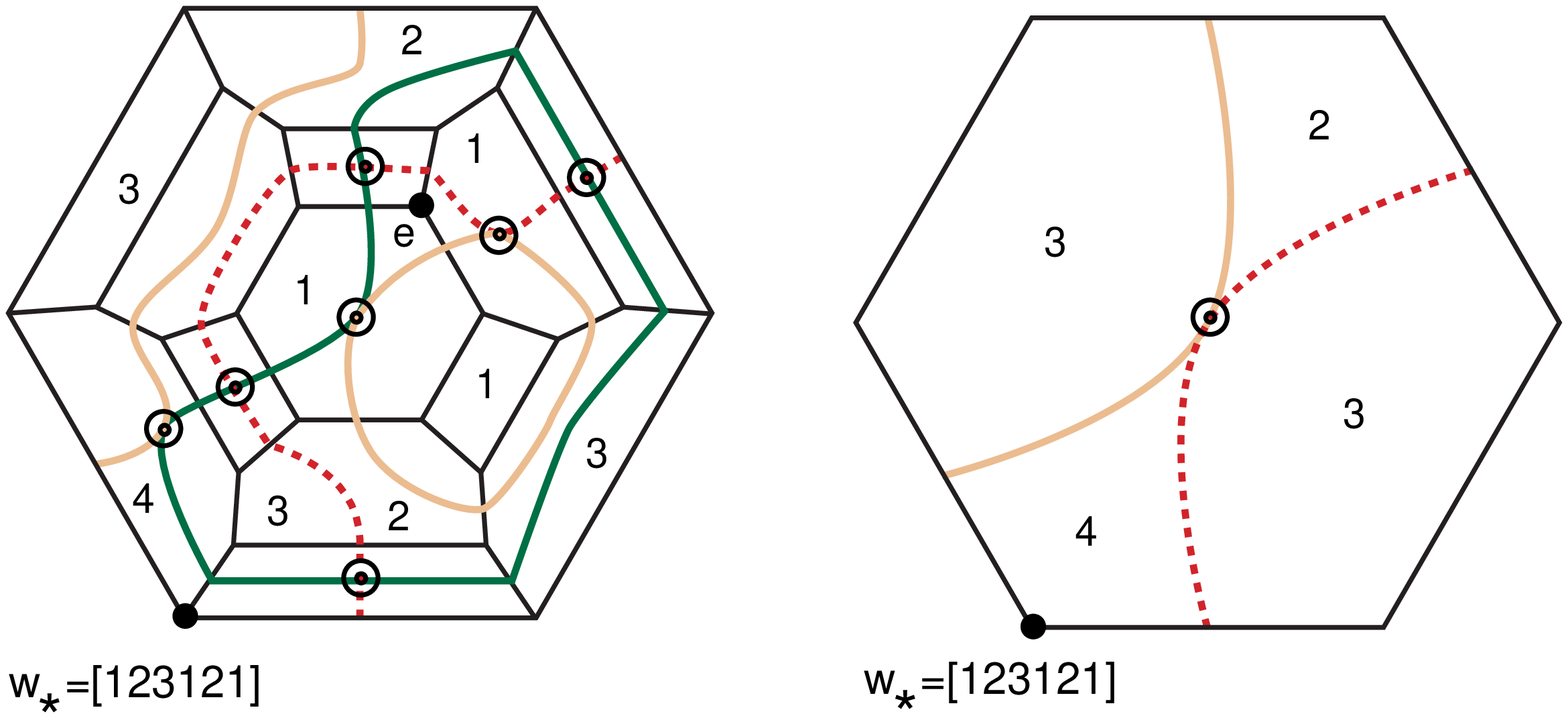}
\caption{The $\Gamma_-$-polytope for type $A_3$ and the Painlev\'e divisors (the right figure
is the back view of the left one). The
Painlev\'e divisors are shown by the dotted curve for ${\mathcal D}_{\{1\}}$, by the light color one for ${\mathcal D}_{\{2\}}$, and by the dark one for ${\mathcal D}_{\{3\}}$. The double circles indicate the divisor
${\mathcal D}_{\{i,j\}}={\mathcal D}_{\{i\}}\cap{\mathcal D}_{\{j\}}$, which are all connected at the center of
the polytope $p_o$. The divisor ${\mathcal D}_{\{2\}}$ has the $A_1$-type singularity
at $p_o$, i.e. a double cone with $t_2^2-t_1t_3=0$.The numbers indicate $\eta(w)$ which are obtained by using any path from $e$ to $w$  along edges of the polytope,  following the direction of the Toda flow, i.e. the Bruhat order.}
\label{A3pol:fig}
\end{figure}

\begin{figure}[t]
\includegraphics[width=9cm]{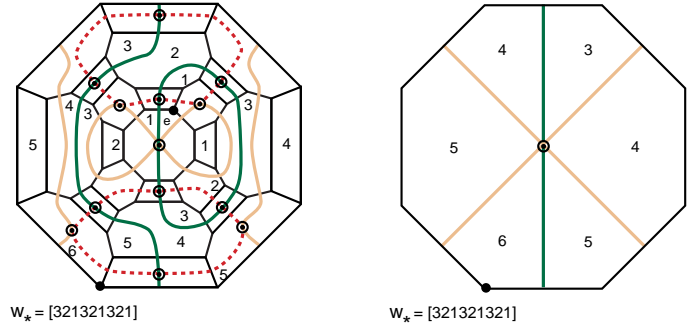}
\caption{The polytope $\Gamma_{-}$ associated to the $B_3$-Toda lattice and the Painlev\'e
divisors. The description of the Painlev\'e divisors are the same as in the case of $A_3$.
The singularities  of ${\mathcal D}_{\{1\}}$ and ${\mathcal D}_{\{3\}}$ are both of $A_1$-type,
while the singularity of ${\mathcal D}_{\{2\}}$ is not isolated, a line singularity attached to two double cones of $A_1$-type (notice two 8-figures of ${\mathcal D}_{\{2\}}$).
The numbers indicate $\eta(w)$. The subsystems on the boundary of $\Gamma_-$ consist of the octagons for $B_2$, the hexagons for $A_2$ and the squares
for $A_1\times A_1$.}
\label{B3pol:fig}
\end{figure}

\begin{figure}[t]
\includegraphics[width=9cm]{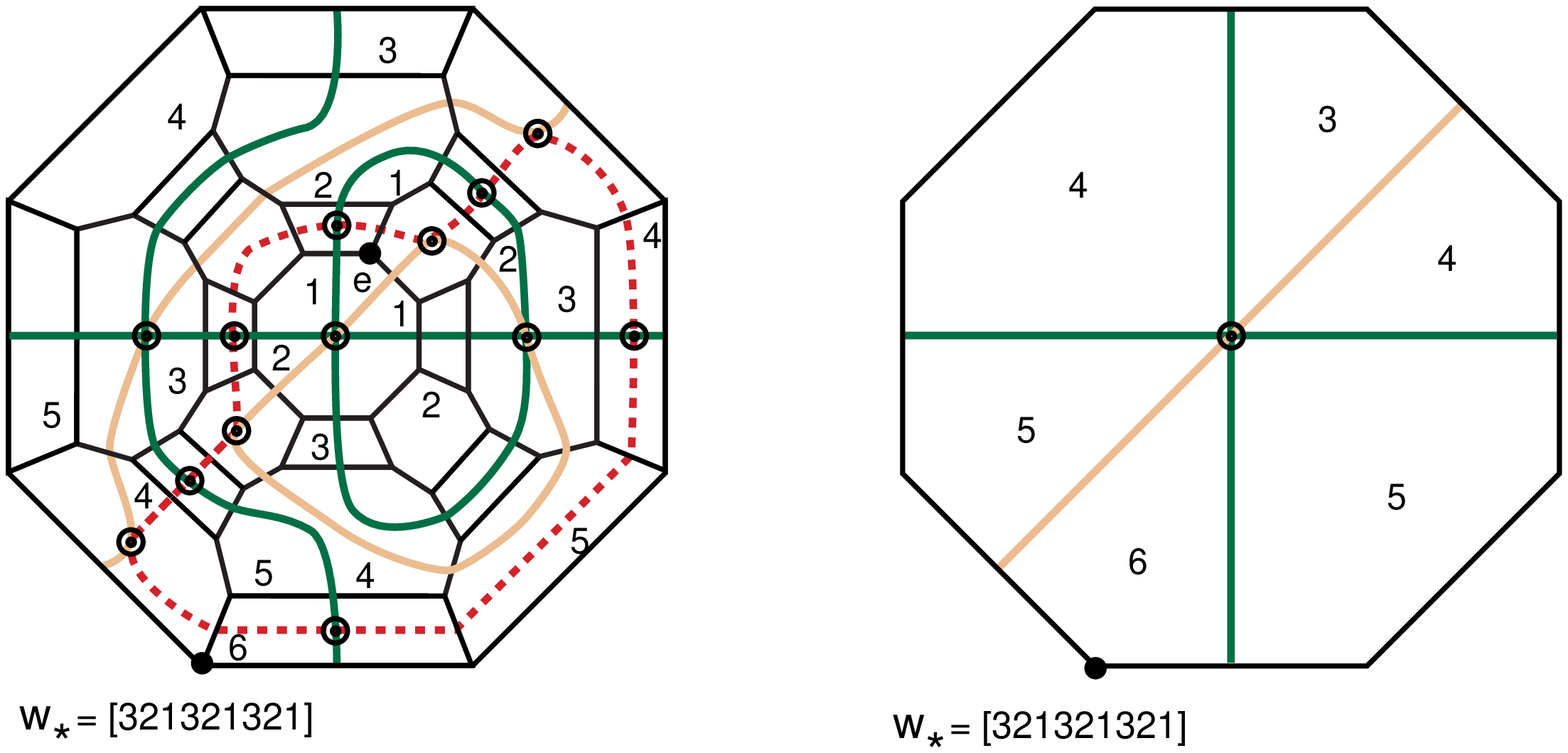}
\caption{The polytope $\Gamma_{-}$ associated to the $C_3$-Toda lattice
and the Painlev\'e divisors. The description of the divisors are again the same as in the case of
$A_3$. Notice that both polytopes for $B_3$ and $C_3$ are the same, but the geometry of the Painlev\'e divisors
are quite different.  The singularity of ${\mathcal D}_{\{2\}}$ is of $A_1$-type, and that of ${\mathcal D}_{\{3\}}$ is a reducible one with $t_1(t_3^2-t_1t_5)=0$.}
\label{C3pol:fig}
\end{figure}

\subsection{Toda lattice: periodic case}\label{periodictoda}
Here we give a brief background of periodic Toda lattice for affine $A_l$ Toda lattice
(the details can be found in \cite{kodama:02}).
The periodic Toda lattice is also give by the Lax equation (\ref{lax}) with
\[
L_P = 
\begin{pmatrix}
b_1 & 1 & 0 & \cdots &0& a_0z^{-1} \\
a_1 & b_2-b_1 & 1 & \cdots &0 & 0 \\
0&a_2&b_3-b_2&\cdots&\cdot&0\\
\vdots & \vdots & \vdots & \ddots & \vdots&\vdots\\
0 & \cdot & \cdot &\cdots & b_{l}-b_{l-1} & 1 \\
z &0 & \cdot &\cdots & a_{l} & -b_l \\
\end{pmatrix}.
\]
The characteristic equation for $L_P$ defines the algebraic curve,
\begin{equation}
\label{characteristicLP}
{\rm det}(L_P-\lambda I)=-\displaystyle{\left(
z+\frac{\prod_{i=0}^l a_i}{ z}-P(\lambda)\right)}=0\,,
\end{equation}
which is a
1-dimensional affine variety on $(\lambda,z) \in {\mathbb C}^2$. Here $P(\lambda)$
is an $(l+1)$-th polynomial of $\lambda$ given by
\begin{equation}
\nonumber
P(\lambda):=\Delta_{l}(\lambda)-a_0\Delta_{l-1}(\lambda),
\end{equation}
where $\Delta_{l}(\lambda)={\rm det}(L-\lambda I)$ for the Lax matrix $L$ of the
nonperiodic Toda lattice,  and $\Delta_{l-1}(\lambda)$ is the determinant $\Delta_l(\lambda)$
after removing the first row and the last column.
Then the polynomial $P(\lambda)$ gives the set of integrals $\{I_k(L_P): k=1,\ldots,l\}$ of the Toda flow, i.e.
\[
\label{P}
P(\lambda)=-(-1)^l\left(\lambda^{l+1}+\sum_{k=1}^l(-1)^kI_k(L_P)\lambda^{l-k}
\right),
\]
There exists an additional integral obtained from the residue with respect to the
spectral parameter $z$, i.e. $I_{0}=\prod_{i=0}^l a_i$. One should note that
the case with $a_0=0$, i.e. $I_{0}=0$ corresponds to the nonperiodic Toda lattice.
Then the isospectral set is defined by
\[
Z_{\mathbb R}^P(\gamma):=\Big\{(a_0,a_1,\cdots,a_l,b_1,\cdots,b_l)\in Z_{\mathbb R}^P:
I_k(L_P)=\gamma_k\in{\mathbb R},~k=0,1,\cdots,l\Big\}
\]
which is the affine part of the compactified manifold ${\hat Z}_{\mathbb
R}^P(\gamma)$ of ${\rm dim}Z^P_{\mathbb R}(\gamma)=l$,
and with a divisor ${\Theta}$ associated to the blow-ups of
$a_k, k=0,1,\cdots,l$, we have \cite{adler:91}, 
\[
\label{ZR}
Z_{\mathbb R}^P(\gamma)={\hat Z}_{\mathbb R}^P(\gamma)\setminus { \mathcal D}_{\{0,1,\ldots,l\}},
\quad {\rm with} \quad {\mathcal D}_{\{0,1,\ldots,l\}}=\bigcup_{k=0}^l\{a_k^{-1}=0\}.
\]
It turns out that the flow of Toda lattice can be described as a trajectory on
the Riemann surface, $y^2=P(\lambda)^2-4I_{0}$, and through the Abel-Jacobi map,
the compactified isospectral manifold ${\hat Z}_{\mathbb R}^P(\gamma)$
can be identified as the real part of the Jacobian
${\mathbb C}^l/\Lambda$ with the lattice $\Lambda$ defined by the period matrix
$\Omega$ associated with the Riemann surface of genus $l$. The divisors
${\mathcal D}_{\{j\}}=\{a_j^{-1}=0\}$ are given by the theta divisor and its translates, that is, the zeros of the Riemann theta function associated with the hyperelliptic Riemann surface, $y^2=P(\lambda)^2-4I_{0}$.
Here we take appropriate values of the integrals $I_k=\gamma_k$ 
(e.g. $I_{0}=\prod_ka_k>0$) so that all the roots of the hyperelliptic
curve are real and distinct , i.e. $y^2=\prod_{k=1}^{2l+2}(\lambda-\lambda_k)$ with $\lambda_k\in{\mathbb R}$. Then the number of connected components in the 
compact manifold $\hat Z^P_{\mathbb R}(\gamma)$ is given by $2^l$, that is,
the real part of the Jacobian consists of $2^l$ number of $l$-dimensional
tori (see \cite{kodama:02}).

In the case of periodic Toda lattice, the signs of $a_k$ are determined by the action
of the {\it affine} Weyl group associated with the affine Kac-Moody algebra (see \cite{kodama:02}).
For example, for the case of $A^{(1)}_l$ Toda lattice with $l\ge 2$, we have
\[
\hat W=\left\langle\
s_0,s_1,\ldots,s_l : \begin{array}{lllll}
s_k^2=e,\quad (s_ks_{k+1})^3=e~{\rm mod}(l+1) \\
(s_ks_j)^2=e,\quad 1<|k-j|<l
\end{array}
\right\rangle
\]
and for $A^{(1)}_1$, we have $
\hat W=\langle s_0,s_1: s_0^2=s_1^2=e\rangle$.
The action of $\hat W$ is defined by the same way as in Definition \ref{act} (see Eq.(4.5) in p.1710 of \cite{kodama:02}), i.e.
for each $s_i\in \hat W$ and $\epsilon_j={\rm sgn}(a_j)$,
\[
s_i: \epsilon_j\longmapsto \epsilon_j\epsilon_i^{-\hat C_{j,i}}
\]
where $\hat C$ is the extended Cartan matrix. 

\begin{figure}[t]
\includegraphics[width=6cm]{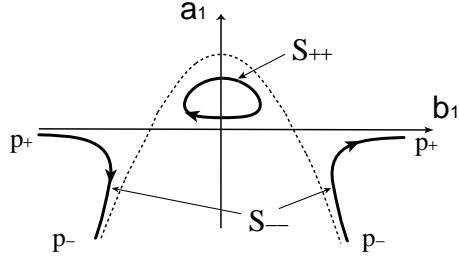}
\caption{The isospectral manifold for the affine $A_1$-Toda lattice.
The solid curves show the invariant curve $I_1=b_1^2+a_1+\frac{\gamma_0}{a_1}=\gamma_1$ for
$\gamma_1>2$ and $\gamma_0=a_0a_1=1$,
and the points $p_{\pm}$ indicate the blow-up points, i.e. the divisors ${\mathcal D}_{\{0\}}$
and ${\mathcal D}_{\{1\}}$. The compactified manifold $\hat Z^P_{\mathbb R}(\gamma)$
is given by a disjoint union of two circles, $S_{++}$ and $S_{--}$, where the signs $(\epsilon_0,\epsilon_1)$
correspond to $\epsilon_k=$sgn$(a_k)$.
The dashed curve indicates the non-periodic limit of the Toda lattice, i.e.
$\gamma_0\to 0$.}
\label{A1-P:fig}
\end{figure}

\begin{Example}
The affine ${\mathfrak{sl}}(2)$ Toda lattice: The Lax matrix is given by
\[
L_P=\left(
\begin{matrix}
b_1 & 1+a_0z^{-1} \\
z+a_1 & -b_1
\end{matrix}
\right)
\]
whose spectral curve defines an elliptic curve,
\begin{equation}
\nonumber
{\rm det}(L_P-\lambda I)=-\displaystyle{\left(z+\frac{a_0a_1}{z}-P(\lambda)\right)}=0.
\end{equation}
The polynomial $P(\lambda)$ is given by
\begin{equation}
\label{g1I2}
P(\lambda)=\lambda^2-I_1, \quad {\rm with} \quad I_1=b_1^2+a_1+
\displaystyle{\frac{1}{a_1}},
\end{equation}
where we have assumed $I_0=a_0a_1=\gamma_0=1$.
Then the affine part of the isospectral manifold is given by
\begin{equation}
\nonumber
Z_{\mathbb R}^P(\gamma)=\Big\{(a_1,b_1)\in {\mathbb R}^2:
b_1^2+a_1+\frac{1}{ a_1}=\gamma_1\in{\mathbb R}\Big\}
\cong \left\{
\begin{matrix}
 S^1\sqcup {\mathbb R}\sqcup {\mathbb R} \quad & {\rm if} ~~ \gamma_1>2 , \\
 {\mathbb R}\sqcup {\mathbb R} \quad & {\rm if} ~~ \gamma_1<2 .
\end{matrix}
\right.
\end{equation}
Two disconnected pieces of $\mathbb R$ are compactified to make a circle $S^1$,
that is, the corresponding solution blows up once in each point $p_+$
or $p_-$, which are the infinite points of the Riemann surface, and total twice in one cycle.
Thus the compactified manifold ${\hat Z}_{\mathbb R}^P(\gamma)$ is just a disjoint
union of two circles, and each circle can be marked by the signs of $a_k$. In Figure \ref{A1-P:fig},
those circles are denoted by $S_{++}$ and $S_{--}$. The flow on the circle $S_{++}$
is complete, and the flow on $S_{--}$ has blow-ups. 
\end{Example}

\subsection{The group $G$ and the group $\check G$}\label{groupG}
We  give here some remarks on the Langlands dual $\check G$ of the real connected Lie group $G$ with ${\rm Lie}(G)=\mathfrak g$ and finite Chevalley groups in connection with our present study:
Let $G({\mathbb C})$ be a connected semisimple Lie group with Lie algebra ${\mathfrak g}_{\mathbb C}={\mathfrak g} + \sqrt{-1}{\mathfrak g}$.
This can be regarded as an algebraic group defined over a finite extension of the field ${\mathbb Q}$  of rational numbers. It is then possible
to consider the group of real points $G({\mathbb R})$. We  denote by $G$ the real connected Lie subgroup of  $G({\mathbb C})$  with Lie algebra  ${\mathfrak g} $. 
Hence $G \subset G({\mathbb R}) \subset  G({\mathbb C})$.  It is also possible to consider this algebraic group
over fields of characteristic $p$, with $p$ a prime. 

We also use the notation $\check G$ to denote any group
associated to $\check {\mathfrak g}$ in the same way as $G$ is associated to ${\mathfrak g}$. 
We will refer to $\check G$ loosely by the term {\it Langlands dual}.  In the case of simple Lie algebras
not of type  $B_l$ or $C_l$ with $l \ge 3$ we may just asume $\check G=G$ in all the statements. 
The Langlands dual will be needed to explain the connection between the blow-ups 
of the Toda lattice associated with ${\mathfrak g}$ and
the cohomology of the real flag manifold for $\check{\mathfrak g}$.

\begin{Example} We consider $G({\mathbb C})=SL(n;{\mathbb C})=\check G ({\mathbb C})$, i.e. the set 
of all the $n\times n$ complex matrices $A$ satisfying
the polynomial equation, ${\rm det}(A)=1$.  Then the complex
solutions of  ${\rm det}(A)=1$ define $G({\mathbb C})$. The group
of real points is, of course, $G({\mathbb R})=SL(n;{\mathbb R})$; and, since this group is connected, $G({\mathbb R})=G$.
There is also an involution $\theta$  given by $\theta(A)=A^*$, the inverse of the
transpose of $A$. We then have  that a maximal compact  Lie subgroup of $G$ is  $K=SO(n)$
 given as the set of matrices satisfying $\theta(A)=A$. 

 In the case of a Lie algebra of type $A_1$ we then have two possibilities for
$G$ namely $G({\mathbb R})=G=SL(2;{\mathbb R})$ or $G({\mathbb R})=G={\rm Ad}SL(2;{\mathbb R})$ the adjoint group. 
This depends on whether we pick $G({\mathbb C})=SL(2;{\mathbb C})$ or  $G({\mathbb C})={\rm Ad}SL(2;{\mathbb C})$.
If  we now let  $G({\mathbb R})=G=SL(2;{\mathbb R})$ then both $\check G=SL(2;{\mathbb R})$  or  $\check G={\rm Ad}SL(2;{\mathbb R})$
are possible.

The equation ${\rm det}(A)=1$ has integral coefficients.  The integral coefficients
 make it possible to reduce modulo a prime $p$. Let $q$ be a power of a prime
$p$ and $k_q$ an algebraic closure of the finite field with $q$ elements denoted  ${\mathbb F}_q$. 
We may then consider  $G(k_q)=SL(n;k_q)$ i.e.
 the solutions in $k_q$ of ${\rm det}(A)=1$.  
By also reducing modulo $p$ the involution $\theta(A)=A^*$
we obtain $SO(n;k_q)$ as the set of 
fixed points and then the finite group  $SO(n;{\mathbb F}_q)$.
In the simplest example of $n=2$, $SO(n; k_q)$  consists of $2\times 2$ matrices
$\begin{pmatrix}
x & y   \\
-y & x \\
\end{pmatrix}$ satisfying $x^2+y^2=1$ with $x,y \in k_q$.
The number of points of the finite group $SO(2;{\mathbb F}_q)$ is then 
the number of  solutions of $x^2+y^2=1$ over the field ${\mathbb F}_q$.
The answer is given by the polynomial $q-1$ if we assume that $\sqrt{-1}\in {\mathbb F}_q$.
\end{Example}

\section{The polynomial $p(q)$ as alternating sum of the blow-ups}\label{pqpolynomial}
In the nonperiodic case, we introduce polynomials in terms of the numbers $\eta(w,\epsilon)$
in Definition \ref{eta}, which plays a key role for counting the number of blow-ups.
For each $\Gamma_{\epsilon}$-polytope we then define $p_\epsilon(q)$.

\begin{Definition} We define a polynomial, the alternating sum of the blow-ups,
\begin{equation}
\label{pq}
 p_\epsilon(q)=(-1)^{l(w_*)}\sum_{w\in W} (-1)^{l(w)}q^{\eta(w,\epsilon)}\,,
\end{equation}
where $w_*$  is the
longest element of the Weyl group and $l(w)$ indicates the length of $w$. When $\epsilon=(-\cdots -)$, we 
simply denote this polynomial by $p(q)$.
\end{Definition}

One can easily show that $p_\epsilon (q)=0$ unless $\epsilon=(-\cdots -)$. Hence the only relevant
polynomial here is $p(q)$, the polynomial for the $\Gamma_{-\cdots -}$-polytope. 
This  corresponds  to the fact that the {\it rational}
cohomology of $K$ and ${\mathcal B}=K/T$  actually agree.
 Recall that we are dealing only with the case
when the Lie algebra is split and the group $T$ is then a finite group.
The polynomial $p_\epsilon (q)$ corresponds to a Lefschetz number
for the Frobenius action over a field of positive characteristic 
for cohomology with local coefficients \cite{casian:05}.  When $q=1$ these polynomials all vanish,
including $p(q)$. This reflects the fact that the Euler characteristic of $K$ is zero.

As we will explain below that the numbers $\eta(w,\epsilon)$ are deeply tied up with
the cohomology of the flag manifold. The polynomial $p(q)$ contains
all the information regarding the cohomology ring of the compact
Lie group $K$. We discuss below our approach to obtaining this strange relations among
the Toda lattice, the flag manifold and $K$.  One of our main results
is the computation of the polynomials $p(q)$ in terms of $K$. However
there is an important technical point here. The polynomial $p(q)$
from the Toda Lattice associated to a Lie algebra ${\mathfrak g}$
is given in terms of the maximal compact subgroup of $\check G$.
More precisely, one must consider the group $\check K$ over ${\mathbb F}_q$
and $p(q)$ is given in terms of the polynomial $|\check K({\mathbb F}_q)|$.
If  the Lie algebra is simple and not of type $B_l$ or $C_l$ with $l\ge 3$
then we can just  write $\check G=G$ and $\check K=K$ in all the statements.

\begin{Example} From Figures \ref{hexagon1:fig}
and \ref{G2:fig}, we note that the numbers $\eta(w,\epsilon)$ are constant
on the connected components in a given $\Gamma_{\epsilon}$-polytope.  The polynomials $p(q)$
that are obtained from (\ref{pq}) by counting $\eta(w)$ are $p(q)=q^2-1$ for type $A_2$, and
$p(q)=(q^2-1)^2$ for $G_2$.  
Note that the multiplicities $d$ are the degrees of these polynomials.
In the case of $G_2$ each divisor contributes $2$ to the multiplicity.
Note in Figure \ref{G2:fig} that  the divisors shown in the polytopes $\Gamma_{\epsilon}$ with
$\epsilon=(+-)$ or $(-+)$ come together in the $\Gamma_{--}$-polytope
 giving a  multiplicity $d=4$ at the center $p_o$ of the $\Gamma_{--}$-polytope.
 
 One can also directly compute the polynomial $p(q)$ from Figures \ref{A3pol:fig}, \ref{B3pol:fig} and \ref{C3pol:fig} for the cases of $A_3, B_3$ and $C_3$. Namely we obtain
 $p(q)=(q^2-1)^2$ for $A_3$, $p(q)=(q-1)(q^2-1)(q^3-1)$ for $B_3$, and $p(q)=(q^2-1)^3$ for $C_3$.
The polynomial obtained from the $B_3$ case can be written as $q^{-3}|\check K({\mathbb F}_q)|$ where $\check K=U(3)$ is the maximal
compact subgroup of  the Langlands dual $Sp(3;{\mathbb R})$, that is, $C_3$.  This serves to illustrate the
fact that the polynomial $p(q)$ is related to the Chevalley group that results from a maximal
compact subgroup of the Langlands dual $\check G$.

Note also  that $A_3$ case gives the same polynomial $p(q)=(q^2-1)^2$ obtained in the $G_2$ case. The reason is that the maximal compact group $K$ is in both cases essentially $SU(2)\times SU(2)$. This will be explained in Theorem \ref{Kpq} which gives the general formulae for
 the polynomials $p(q)$ for all the cases.
 \end{Example}

In the affine (periodic)  cases we just consider  $p_\epsilon (q)=\sum_{w\in\hat W} (-1)^{l(w)}q^{\eta(w,\epsilon)}$, which is now given by a power series of $q$.
The numbers $\eta(w,\epsilon )$ are defined similarly but the elements $w\in \hat W$ do not
represents points in the compactified isospectral manifold.  In the concrete examples given here
the universal cover (${\mathbb R}^{l}$)  is subdivided into regions by the divisors, which
are further subdivided so that they are labeled by Weyl group elements
(see Figure \ref{A2aff:fig} below).  The numbers $\eta(w,\epsilon)$
are assigned uniquely to the different regions by counting blow-up points along the Toda trajectories ignoring the direction of the flow. If  a path of the Toda lattice goes  from a region labeled $w$ to
a region labeled $w^\prime$, with $w \le w'$ in the Bruhat order,  and the path crosses $k$ blow-up points, 
then $\eta(w^\prime,\epsilon)=\eta(w,\epsilon)+k$. Setting $\eta(w,\epsilon)=0$ determines all the other numbers uniquely. Note that  there may not be a concrete path  going from the region labeled $e$ to the region labeled $w$ but, still,  the number  $\eta(w,\epsilon)$ is still being  determined by 
counting blow-ups along trajectories of the Toda lattice.

\section{A graph associated to the blow-ups of the Toda lattice}

The following graph ${\mathcal G}_\epsilon$ was originally  motivated by the problem of computing
the number of connected components in the $\Gamma_{\epsilon}$-polytope. This problem
is analogous to the problem of computing the intersection of two opposite 
top dimensional Bruhat cells in the case of a real flag manifold (e.g. see \cite{shapiro:97,rietsch:97,
zelevinsky:00}).
We then observed that in all examples
this was the graph of incidence numbers for a real flag manifold (see \cite{casian:99}).
In the affine cases we may consider
$\hat G=G({\mathbb R}[t, t^{-1}])$ instead of
 $G$ and $\hat B=\{ f\in G({\mathbb R}[t, t^{-1}]) : f(0) \in B \}$ instead of $B$.
Thus $\hat {\mathcal B}=\hat G/ \hat B$.

\begin{Definition}
\label{graph}
For a fixed $\epsilon=(\epsilon_1\ldots \epsilon_l)$, we associate a graph
${\mathcal G}_\epsilon$ to the blow-ups of the
Toda lattice. The graph consists of vertices labeled by the elements
 of the Weyl group $W$, i.e. the vertices of the $\Gamma_{\epsilon}$-polytope, and 
 oriented edges $\Rightarrow$. The edges are defined as follows:
For any $w_1 , w_2 \in W$, there exists an edge between $w_1$ and $w_2$,
\[
w_1\Rightarrow w_2\quad {\rm iff}\quad
\left\{\begin{array}{llll}
{\rm a)}~ w_1\le w_2\,\,({\rm Bruhat~order})\\
{\rm b)}~l(w_2)=l(w_1)+1\,, \\
{\rm c)}~\eta(w_1,\epsilon)=\eta(w_2,\epsilon)\,,\\
{\rm d)}~w_1^{-1}\epsilon=w_2^{-1}\epsilon\,.
\end{array}\right.
\] 
When $\epsilon=(-\ldots -)$, we simply denote ${\mathcal G}={\mathcal G}_\epsilon$.
This graph also makes sense in the periodic case.
\end{Definition}

Note that in non-periodic cases $w_1\Rightarrow w_2$ implies that in the polytope there is a path from $w_1$-vertex to
$w_2$-vertex without crossing a Painlev\'e divisor (no blow-up). Namely, in this case,
$w_1\Rightarrow w_2$  means that the vertices associated to $w_1$
and $w_2$ belong to the {\it same} connected component of the hexagon (when blow-ups
are removed). 

In many cases the graph ${\mathcal G}_\epsilon$  accomplishes the job of joining together in its
connected components exactly  those vertices $w$ of the polytope $\Gamma_\epsilon$
belonging to the same connected components.  This will be the case in the example considered below. 

\begin{Example} \label{exam}
In the case of $A_2$, we have $s_1(--)=(-+),\, s_2(-+)=(-+)$ which implies
$(--) \to (-+)\, {\overset{s_2}{\Rightarrow} }\,(-+)$
and  $\eta(s_1)=1$, $\eta(s_1s_2)=1$. Therefore the graph ${\mathcal G}$
which encodes blow-up information in  $\Gamma_{--}$ of  Figure 2  is ($s_i $ is replaced with $i$):


$$\begin{matrix} 
{}                     &{e}& {}               & {} & {}& : & q^0 \\
1                     &{} & 2                  & {} &{}& : & q^1 \\
 \Downarrow&{}&\Downarrow & {} & {} &{}& \\
  12                &{}& 21                  & {} &{}& : & q^1\\
   {}&{121}&{ }                      &{}&{}& : & q^2\\
    \end{matrix} $$

Here we have also  listed the monomials  $q^{\eta (w)}$ (in the variable $q$) associated to 
 representatives of the integral cohomology ($w\to \eta (w) \to q^{\eta (w)}$).
As already noted, the vertices of the hexagon $\Gamma_{--}$ belonging
to a connected component in Figure \ref {hexagon1:fig} form a connected
component of this graph (see also Figure \ref{A1:fig}). This graph classifying 
connected components in the hexagon minus the blow-ups, agrees
with the graph in p. 465 of  \cite{casian:99} which is defined very differently  in terms of
incidence numbers. The graph of incidence numbers 
 gives rise to a chain complex by replacing the edges $\Rightarrow$ with
multiplication by $2$,
$${\mathbb Q}\langle e\rangle~ {\overset{\delta_o}{\longrightarrow} }~
{\mathbb Q}\langle s_1\rangle \oplus {\mathbb Q}\langle s_2\rangle ~{\overset{\delta_1}{\longrightarrow}}~ {\mathbb Q}\langle s_1s_2\rangle\oplus {\mathbb Q}\langle s_2s_1\rangle~ {\overset{\delta_2}{\longrightarrow}} ~{\mathbb Q}\langle s_1s_2s_1\rangle\,.  $$
Here $\langle w\rangle$ is the Bruhat cell associated to the element $w\in W$.
The only non-zero map is $\delta_1$ given by a diagonal matrix  with $2$'s in the diagonal corresponding to the $\Rightarrow$. The rational cohomology that results is
$$\left\{ \begin{matrix}
H^0(G/B, \mathbb Q )&= &\mathbb Q &  {} &: &{q^0}\\  
H^1(G/B, \mathbb Q )&=& 0                & {} &:& {q^1}\\
H^2(G/B,\mathbb Q )&=&0& {} &:& {q^1}\\
H^3(G/B,\mathbb Q)&=&\mathbb Q & {} &:&  {q^2}\\
\end{matrix}\right. $$

Over the rationals this just gives the cohomology of $K=SO(3)$. The alternating sum of the 
$q^{\eta (w)}$ produces the polynomial $p(q)=q^2-1$. Now $p(q)$  multiplied by $q^{r}$ with
$r={\rm dim}(K)-{\rm deg}(p(q))=1$ gives the number of points
of $SO(3)$ over a field with $q$ elements  ($q$ is a power of an odd prime $p$).
The explanation for this is that the $q^{\eta (w) }$ listed
can be shown to be  Frobenius eigenvalues in etale cohomology
of the appropriate varieties reduced to a field of positive characteristic.
When this is taken into account we see that, in fact, more than
the cohomology of $K/T$ or $K$, we are obtaining  etale $\bar {\mathbb Q}_l$
cohomology over a field of  positive characteristic, including
Frobenius eigenvalues. All derived from the structure of the Toda
lattice and its blow-up points.

\end{Example}

If we start with $\epsilon=(-+)$ then we obtain the edges
$e\Rightarrow s_2$, $s_1\Rightarrow s_2s_1$, $s_1s_2\Rightarrow s_1s_2s_1$.
This time the vertices of the hexagon $\Gamma_{-+}$ belonging
to a connected component in Figure \ref {hexagon1:fig} form a connected
commponent of this graph.
The graph for $\Gamma_{-+}$ now corresponds to the graph of incidence numbers computing
cohomology with local coefficients. The local system  ${\mathcal L}$ can be described 
by the signs $(-+)$. The $-$ sign indicates that along a circle in $G/B$
that corresponds to $s_1$ the local system is constant, and the second $+$
that along a circle corresponding to $s_2$ it is non-trivial. 
With these local coefficients $H^*(G/B;{\mathcal L})=0$,
which implies $p_{-+}(q)=0$ (see Section \ref{pqpolynomial}).
Similarly for $\Gamma_{+-}$ and $\Gamma_{++}$. 

\begin{figure}[t!]
\includegraphics[width=10.5cm]{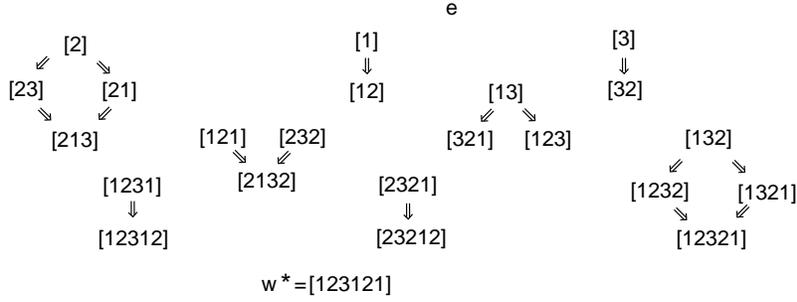}
\caption{The graph $\mathcal G$ of the real flag manifold for type $A_3$.
The Bruhat cells $NwB/B$ are denoted by $[ij\ldots k]$ for $w=s_is_j\ldots s_k$.
The incidence numbers associated with the edges $\Rightarrow$ are $\pm 2$ (see also
Example (8.1) in \cite{casian:99}).
There are 10 connected components in this graph corresponding to the 10 connected components in $\Gamma_{-}$ after blow-up points are removed. 
This graph can be also obtained from Figure \ref{A3pol:fig}.}
\label{A3inc:fig}
\end{figure}

 In the case of a Lie algebra of type $A_3$ we obtain the graph $\mathcal G$ in Figure \ref{A3inc:fig}.
This graph corresponds to the  polytope in Figure \ref{A3pol:fig}  as separated into connected components by the divisors shown.  To determine the number $\eta(w)$ for any given $w$ it is enough to go from $e$ to $w$ along any path along the boundary   corresponding to a reduced expression $w=s_{n_1}\cdots s_{n_r}$ and count the number of intersections with the divisors. In Figure \ref{A3pol:fig}, we show the path
following the expression $w_*=[123121]$, i.e. $e\to s_1\to s_1s_2\to s_1s_2s_3\to \cdots \to w_*$, with the arrows on the edges of the polytope.

We then state the following theorem showing the equivalence between the connected components
in the polytopes $\Gamma_{\epsilon}$ and the graphs ${\mathcal G}_{\epsilon}$ of the
incidence numbers defined in \cite{casian:99}. A proof of the theorem can be found in \cite{casian:05}.
We state the case of  ${\mathcal G}$ but the Remark below explains the general version of the theorem.

\begin{Theorem}\label{mainT} The graph ${\mathcal G}$ is the graph of incidence numbers
for the  cohomology of  the real flag manifold $\check {\mathcal B}$  in terms of 
the Bruhat cells. \end{Theorem}

\begin{Remark} Note that while the compactified  isospectral  manifold of the Toda lattice lives inside $ {\mathcal B}$ through the companion embedding \cite{casian:02},
this theorem mysteriously  involves the  topology of the flag manifold associated to the Langlands dual, namely   $\check {\mathcal B}$.  In Section \ref{singularstructure} we give an explicit evidence of this duality in the computation of the multiplicity of the $\tau$-function at $p_o$.
\end{Remark}
 
 \begin{Remark}  In general each $\check K$-equivariant local system  ${\mathcal L}$ on $\check {\mathcal B}$ corresponds to an $\epsilon$ and  ${\mathcal G}_\epsilon$ is the graph of incidence numbers for cohomology of $\check {\mathcal B}$ with twisted coefficients in ${\mathcal L}$. There may be some signs $\epsilon$ such that no $\check K$-equivariant local system ${\mathcal L}$ corresponds to it.  For instance  if $G=SL(4;{\mathbb R})$ this is the case. Moreover, Theorem \ref{mainT} can be proved for integral
coefficients.
  \end{Remark}

\begin{figure}[t]
\includegraphics[width=9cm]{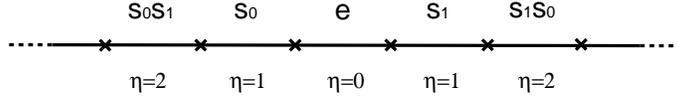}
\caption{The compactified isospectral manifold for the periodic Toda Lattice in the case of $A^{(1)}_1$ is a union of two circles $S_{++}$ and $S_{--}$ in Figure \ref{A1-P:fig}.
This is the universal cover of $S_{--}$. The $\times$ are the divisors and we include the number of blow-ups $\eta(w)$. Each line segment separated by the blow-ups is marked by a unique element
of  the affine Weyl group, $\hat W=\langle s_0,s_1\rangle$. The graph ${\mathcal G}_{--}$ then consists
of disconnected vertices marked by the lements of $\hat W$, and the function $p(q)$
is given by $p(q)=1-2q+2q^2-\ldots=\frac{1-q}{1+q}$.}
\label{A1aff:fig}
\end{figure}

\begin{figure}[t]
\includegraphics[width=11cm]{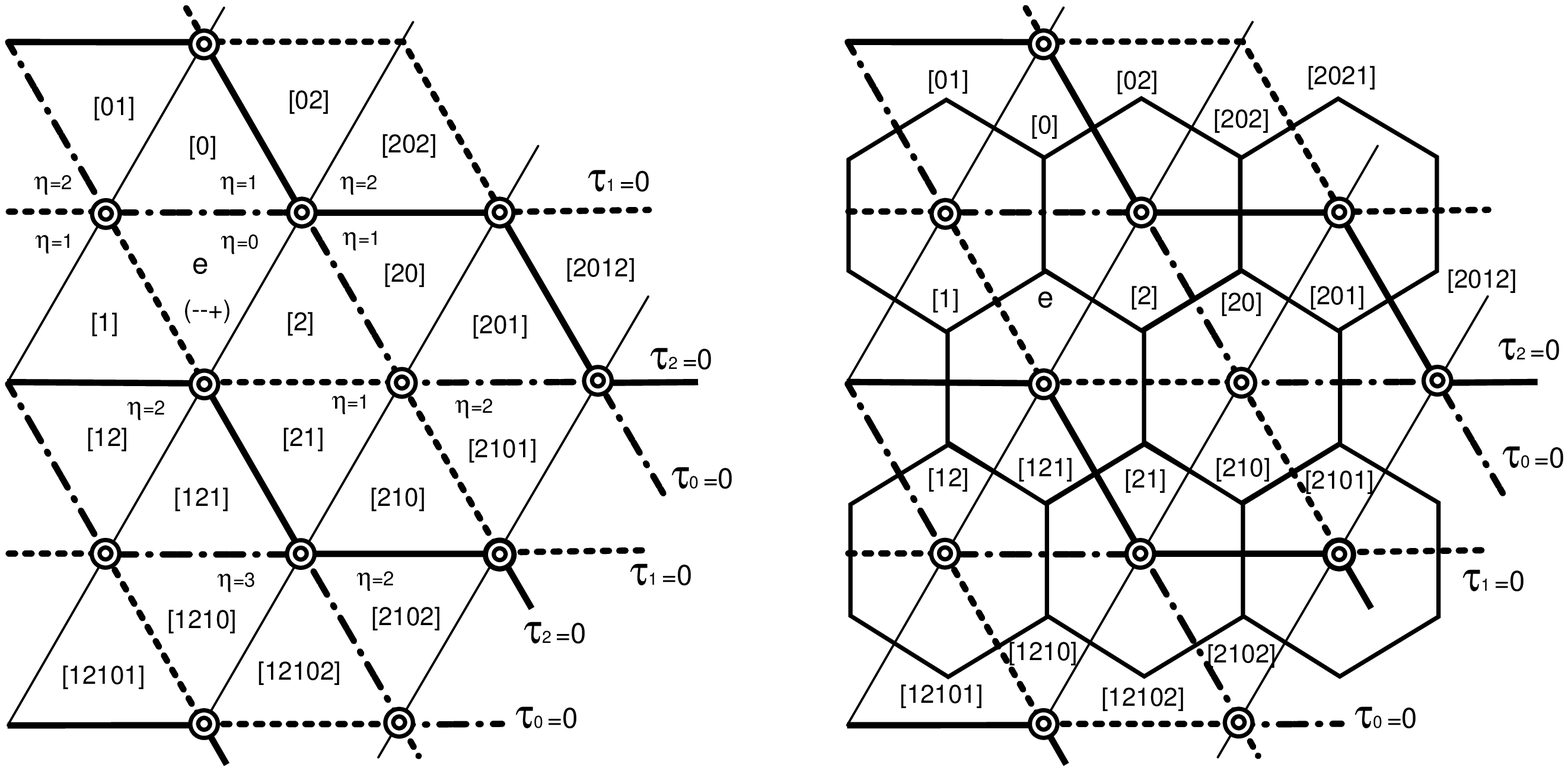}
\caption{The universal cover of one of the tori of the real part of the Jacobian for $A^{(1)}_2$ Toda lattice. Each triangular region is marked by the element
$w$ of the affine Weyl group, e.g. $[02]=s_0s_2$. The right figure shows that each element of
the affine Weyl group is identified as a unique vertex of the honeycomb of the dual graph of
the left figure. Each hexagon represents the Weyl group for $A_2$. Then some of the chamber walls correspond to the divisors (compare with the hexagons in Figure \ref{hexagon1:fig}), and they are indicated here with a solid curve ($\tau_2=0$)
or a dotted curve ($\tau_1=0$) or a dot-dashed curve ($\tau_0=0$).  On a neighborhood of
each double circle ($\tau_k=\tau_{k+1}=0$, $k$ (mod $3$)), the $\tau_k$-functions can be expressed
as the Schur polynomials of $A_2$-nilpotent Toda lattice.
We also include the number of blow-ups $\eta(w)$ starting from $\eta(e)=0$.}
\label{A2aff:fig}
\end{figure}

We now discuss the periodic cases through a couple of examples.
In the periodic cases, the Toda flow can be described as a flow on a Riemann surface, 
and the real solutions are determined by the real part of the surface (see Section \ref{periodictoda}). 
The compactified isospectral manifold of the periodic Toda Lattice consists of union of tori,  the real part of the Jacobian (see \cite{audin:94, kodama:02}).
 The theta divisors are the blow-ups of the Toda Lattice.  For each torus we consider its universal cover
and we need to  pull-back the Toda flow to this universal cover before counting blow-up points.  

\begin{Example}  In the periodic
Toda lattice associated to affine Kac-Moody algebra $A^{(1)}_1$, we have a compactified isospectral manifold consisting of two disjoint circles, i.e. $S_{++}$ and $S_{--}$ in Figure \ref{A1-P:fig}.
Since the only possible boundaries in the chain complex that computes integral cohomology
involve $w$ and $ws_i$ it is easy to compute the graph of incidence numbers. 

In the case of constant coefficients,  the graph of incidence numbers
consists of the affine Weyl group and there are no edges.  The graph
 ${\mathcal G}_{--}$ can easily be computed by counting the number of
blow-ups along the Toda trajectories. We have $\eta(w)=l(w)$ for any $w$ in the
affine Weyl group.  Therefore there are no edges $\Rightarrow$. We obtain
that  the graph  of incidence numbers of the infinite dimensional real flag manifold $\hat G/ \hat B$  that corresponds 
agrees with the graph ${\mathcal G}_{--}$. The case of  $\epsilon=(++)$ corresponds to
a graph of incidence numbers for integral cohomology with twisted coefficients.
Then we have $p(q)= \frac{1-q}{1+q}$ and $p_{++}(q)=0$. We expect that the rational function $p(q)$
contains some information of the (rational) cohomology of the circle $S_{--}$, a real part of the Jacobian, e.g. $H^0={\mathbb Q}$ and $H^1={\mathbb Q}$ by looking at the degrees of
polynomials appearing in $p(q)=\frac{p_1(q)}{p_0(q)}$ (in analogy of the Weil conjecture).
\end{Example}

The case of  $A^{(1)}_2$ has the complication that at least one sign $\epsilon_i$ must be $+$ (recall $I_0=a_0a_1a_2>0$).
Figure \ref{A2aff:fig}  corresponds to the universal cover of a one of these tori
(compare Figure \ref{A2aff:fig} with the Figure 2 in p.1516 of
\cite{audin:94}). One then necessarily  obtains cohomology with twisted coefficients and
 $p_\epsilon (q)=0$.  In the $A_3$ case,  $\epsilon=(----)$ should give rise to integral cohomology. 

We have the following conjecture which clarifies the relation between the integral cohomology
of the real flag manifold and the blow-up structure of the Toda lattice:

\begin{Conjecture} The graph ${\mathcal G}$ is the graph of incidence numbers
for the integral cohomology of the real flag manifold ${\mathcal B}$ ($\hat {\mathcal B}$)  in terms of 
the Bruhat cells.  In general, each $\epsilon$ corresponds to a local system ${\mathcal L}$
in the real flag manifold and the graph ${\mathcal G}_\epsilon$ is the graph
of incidence numbers in the computation of integral cohomology with coefficients in ${\mathcal L}$.
\end{Conjecture}

\begin{Remark} This is an extension of Theorem \ref{mainT} for the finite dimensional case (non-periodic Toda case) which has been completed as Theorem 3.5 in \cite{casian:05}. 
The infinite dimensional case (periodic Toda case)  has several steps which were omitted above. We consider  a real split semisimple Lie algebra ${\mathfrak g}$ and the corresponding
  affine Lie algebra $\hat {\mathfrak g}$.  We also have the Weyl group $W$  and the affine Weyl group $\hat W$. 

First  one considers the universal cover of the compactified isospectral manifold. Since it consists of several tori of dimension $l$, we end up 
with several copies of  $l$ dimensional vector space which we could identify with the Cartan subalgebra ${\mathfrak h}$ . The idea here will be that
the structure of these universal covers  (of each tori), and the pull-back of the Painlev\'e divisors, resemble ${\mathfrak h}$ together with certain walls of an
action of  $\hat W$ on ${\mathfrak h}$.  Step 1 is a division of each of the universal covers into regions parametrized by $\hat W$.  This is a refinement of the
division into connected components by the Painlev\'e divisors pulled back to the universal cover. Moreover, the region 
assigned to the identity has a sign $\epsilon$ and the region corresponding to $w$ will have a sign obtained by applying $w$ to the sign $\epsilon$.
Step 2 is that $\eta (w)$, the number of blow-ups associated to each region is well defined.  This now allows us to define the graph ${\mathcal G}_\epsilon$.
\end{Remark}

\section{The real flag manifold $K/T$ and the variety ${\mathcal O}_o$}

For each real split simple Lie algebra ${\mathfrak g}$  we can consider a connected Lie group $G$ and a maximal compact Lie subgroup $K$ which we 
assume is the set of fixed point sets of a Cartan involution $\theta$. Moreover  in  the appropriate context of algebraic groups, all these objects can also be considered over a field $k_q$, an algebraic closure of a finite field ${\mathbb F}_q$ with $q$ elements.  We just give a simple example and recall  that the real flag manifold can be replaced with a complex manifold ${\mathcal O}_o$ which has the same
homotopy type. This is just the unique open dense   $K({\mathbb C})$ orbit in  the complex flag manifold ${\mathcal B}_{\mathbb C}$.  
 The advantage of this is that we will be able to consider  an analogue of the
real flag manifold that makes sense over fields of positive characteristic.  Note that the complex flag manifold
$G({\mathbb C})/ B({\mathbb C})$ cannot be chosen to play this role since the topologies  of $G({\mathbb C})/ B({\mathbb C})$ and $G/B$ are
very different. For instance, in the $A_1$ case $G/B$ is a circle, i.e. ${\mathbb P}^1$, and $G({\mathbb C})/ B({\mathbb C})$ is a two dimensional sphere, i.e. ${\mathbb{CP}}^1$.

\begin{Example}\label{suoneone} Consider the group ${\bf}= SL(2; k)$.  We consider
the Cartan involution  $\theta$ given
as follows: 
 $\theta(g)= 
\begin{pmatrix}
1 & 0   \\
0 & -1 \\
\end{pmatrix} g 
\begin{pmatrix}
1 & 0   \\
0 & -1 \\
\end{pmatrix}$. Therefore $K(k)$ consists of the diagonal matrices  $\left\{ {\rm diag}(z, z^{-1}) : z\in k \right\}$. The choice of $K$ corresponds to considering
 the subgroup $SU(1,1)$ of $SL(2; {\mathbb C})$ rather than
$SL(2; {\mathbb R})$.  The flag manifold is ${\mathbb P}^1 = k \cup \{ \infty \}$ and the action of 
$K(k)$ is $g(z)\cdot y = z^2y$ for $g(z)\in K(k)$ with $z\in k^*$. The $K(k)$ orbits are $k^*$, $\{ 0 \}$ and $\{ \infty \}$.
Hence for $k={\mathbb C}$, we have ${\mathcal O }_o={\mathbb C}^*$. This orbit has the homotopy type of 
the real flag manifold of $SU(1,1)$, namely a circle $S^1$.
\end{Example}

\section{Computation of $p(q)$ in terms of finite Chevalley groups} 

The Chevalley groups that will be related to the Toda Lattice are the finite groups
$ \check K({\mathbb F}_q)$.   Since Langlands duality  $ {\mathfrak g} \to \check  {\mathfrak g}$ only affects type $B$ and $C$
we will,  for the sake of notational simplicity,  ignore this technical detail, and assume that our Lie algebra
does not contain factors of type $B$ or $C$. 

Note that so far we have explained a relation of the Toda lattice with the cohomology
of $G/B$ and not  with the cohomology of $ K$. However, it turns out that 
$H^*(G/B; {\mathbb Q})=H^*(K;{\mathbb Q})$
(Proposition 6.3 of   \cite{casian:05}).  This is equivalent to showing that for $K$-equivariant
twisted coefficients  ${\mathcal L}$,  $H^*(G/B; {\mathcal L})=0$. This, it turns out, is the reason that
$p_\epsilon (q)=0$ whenever $\epsilon$ contains at least one $+$ sign.

The main idea now is to rely on a Lefschetz Fixed Point Theorem.  
The polynomial $p(q)$ is defined as an alternating sum of powers $q^{\eta (w)}$  the polynomial that computes
the order of the finite group $K({\mathbb F}_q)$ is also given as a polynomial in  $q$.
This second statement can be seen by direct computation,  but it can also be seen as
a consequence of a Lefschetz Fixed Point Theorem. This second way of looking at it is
much more complicated than the direct computation  but it will be more useful for us in this case. This
 requires that we first  consider  the  group $K(k_q)$ with $k_q$ an algebraic closure of the finite
field ${\mathbb F}_q$. Then we need to consider its cohomology, but, since  this algebraic group is not
over ${\mathbb C}$ or ${\mathbb R}$,  ordinary cohomology does not
make sense. We then  use  {\it  etale cohomology}. In etale cohomology
there is also an action induced by the  Frobenius map. The  eigenvalues resulting from this action are, in this case, powers of $q$.
Using a version of the Lefschetz Fixed Point Theorem applied to the Frobenius map $Fr$
one obtains that  $|K({\mathbb F}_q)|$ is an alternating sum of powers of $q$.
Since the  Lefschetz Fixed Point Theorem involves not cohomology, but cohomology
with compact supports,  we will not get exactly  $|K({\mathbb F}_q)|$ but rather 
$q^{-r}|K({\mathbb F}_q)|$ for some $r$. The main point will then be that these
Frobenius eigenvalues are exactly the numbers $q^{\eta (w)}$
obtained from the Toda Lattice.

\begin{Example} 

In the case of $K=SO(2)$ i.e. of a circle, we can give a characteristic zero analogue of this explanation
involving the Frobenius map:  Let us consider the map,
$\Phi_q:S^1 \to S^1$, $z\mapsto z^q$, a map of degree $q$.  Then we have
\[
\left\{\begin{array}{llll}
H_0(S^1;{\mathbb Q})={\mathbb Q}\quad &:\quad q^0\\ 
H_1( S^1;{\mathbb Q})={\mathbb Q}\quad &:\quad q^1
\end{array}\right.
\]
and now the number of fixed points is $q-1$, i.e. the number of non-zero roots of $z^q=z$, or as in
the Lefschetz fixed point theorem,
\[
L(\Phi_q)={\rm Tr}\left((\Phi_q)_*|_{ H_1( S^1;{\mathbb Q})}\right)- {\rm Tr}\left((\Phi_q)_*|_{H_0(S^1;{\mathbb Q})}\right)=q-1\,.
\]

If we replace $SO(2)$ with its complexification, we obtain ${\mathbb C}^*={\mathbb C}\setminus \{ 0 \}$. 
Reduction to positive characteristic means we consider $k_q^*=k_q \setminus \{ 0 \}$.  Now the Frobenius map
is $Fr(z)=z^q$ and the fixed points of this map are the ${\mathbb F}_q$ points i.e., the elements in
${\mathbb F}_q\setminus \{ 0 \}$. We have then  $q-1=|K({\mathbb F}_q)|$. The polynomial $q-1$ is
just $p(q)$ obtained in the case of $A_1$ from the Toda lattice. The number of blow-ups  in Figure \ref {A1:fig}
for $\Gamma_{-}$ is $1$ and $q^1$ is the Frobenius eigenvalue  appearing in etale
cohomology with proper supports of $k_q^*$,  corresponding to  the degree of the map $\Phi_q$.

\end{Example}

We now expand a little  on this explanation 
and the full details are in  \cite{casian:05}: Recall that there is  filtration by Bruhat cells,
\[ 
\emptyset\, \subset \, {\mathcal B}_0\, \subset  \, {\mathcal B}_{1} \, \subset \, \cdots\, \subset \, {\mathcal B}_{l(w_*)}= {G/B}
\]
where ${\mathcal B}_j:=\cup_{l(w)\le j}{NwB/B}$.
There is a similar filtration of ${\mathcal Y}_o$:  $\cdots \subset {\mathcal Y}_{j} \subset {\mathcal Y}_{j+1}\subset \cdots $ given by intersection of ${\mathcal O}_o$
 with $N({\mathbb C})$ cells
inside $G({\mathbb C})/B({\mathbb C})$.  We obtain  coboundary maps:
$H^j({\mathcal Y}_j, {\mathcal Y}_{j-1}; {\mathbb C}) \to 
H^{j+1}({\mathcal Y}_{j+1}, {\mathcal Y}_{j}; {\mathbb C})$ which
give rise to a chain complex computing the cohomology of ${\mathcal O}_o$.
For example in the case of $SU(1,1)$ above,  ${\mathcal Y}_0=\{ \infty \}$ and ${\mathcal Y}_1={\mathbb C} \setminus \{ 0 \}$. 
Each $w$ corresponds to a dual of a Bruhat cell and  contributes to $H^{l(w)}({\mathcal Y}_{l(w)}, {\mathcal Y}_{l(w)-1}; {\mathbb C})$
giving rise to a cohomology class $[w]_{\mathbb C}$.
This can be done with etale cohomology with coefficients in $\bar {\mathbb Q}_m$
and a field of positive characteristic  where $m$ is relatively prime to $p$ and $p\not=2$
and $x^2+1$ factors over ${\mathbb F}_q$.  In this case a Frobenius action arises
in cohomology.

We have the following proposition  for the Frobenius eigenvalue of the cohomology
class $[w]_{k_q}$ in etale cohomology over a field $k_q$ algebraic closure of ${\mathbb F}_q$ of characteristic $p$. This assumes $\bar {\mathbb Q}_m$ coefficients.

\begin{Proposition}\label{frob}
The cohomology class $[w]_{k_q}$ in $H^{l(w)} ({\mathcal Y}_{l(w)}, {\mathcal Y}_{l(w)-1}; \bar {\mathbb Q}_m)$ corresponding to $w\in W$
has Frobenius eigenvalue given by  $q^{\eta(w)}$. 
\end{Proposition}

If one really tries to see where this comes from, this really is Corollary 5.1 in   \cite{casian:05}. The computation of $\eta (w)$
agrees with the computation of certain coefficients that arise from a module over the Hecke algebra introduced
by Lusztig and Vogan in \cite{lusztig:83} .   These coefficients are already  Frobenius eigenvalues in that paper
and are seen to come from the cohomology of the flag manifold and thus of $K(k_q)$.

\begin{Example} In the case of $SU(1,1)$ which was done earlier, we obtain
Frobenius eigenvalues $1$ and $q$ respectively in the degree $0$ and $1$.
\end{Example}

This proposition is obtained by expressing Proposition 9.5  in \cite{casian:99} in terms
of new notation motivated by the Toda lattice. 
The upshot of this is that  the number of blow-up associated to a vertex $w$ in
the Toda lattice and the Frobenius eigenvalue of
 $[w]_{k}$ in $H^{l(w)}({\mathcal Y}(k)_{l(w)}, {\mathcal Y}(k)_{l(w)-1}; \bar {\mathbb Q}_m)$
 are given by the same formula in terms of $\eta (w)$. 
  Then applying the Lefschetz Fixed Point  Theorem to the 
 Frobenius map to count
 the number of ${\mathbb F}_q$ points.  We obtain the following Theorem
 (Theorem 6.5 in \cite{casian:05}):
 
 \begin{Theorem} 
\label{Kpq}
  The polynomial $p(q)$ satisfies  $ p(q)=q^{-r}| \check K({\mathbb F}_q)|$ with
  $r={\rm dim} (\check K)- {\rm deg} (p(q))$. Moreover $p(q)$ factors as   $p(q)= \Pi_{i=1}^g (q^{d_i}-1)$
where $g$ is the rank of $\check K$.  The polynomial $p(q)$  is given by the following explicit formulas:
 
  $A_{l}$:  $\check K=SO(l+1)$,
 \begin{itemize}
 \item[] $l$ even: $p(q)=(q^2-1)(q^4-1)\cdots (q^{l-2}-1)(q^{l}-1),\quad \quad g=l/2$ 
 \item[] $l$ odd:  $p(q)=(q^2-1)(q^4-1) \cdots (q^{l-3}-1)(q^{l-1}-1)(q^g-1),\quad g=(l+1)/2$ 
 \end{itemize}

 $B_l$:  $\check K=U(l)$,
 \begin{itemize}
\item[]  $p(q)=(q-1)(q^2-1)(q^3-1) \cdots (q^l-1),\quad g=l$ \par
\end{itemize}

$C_l$:  $\check K=SO(l)\times SO(l+1)$,
 \begin{itemize}
\item[]  $l$ even:  $p(q)=(q^2-1)^2(q^4-1)^2 \cdots (q^{l-2}-1)^2(q^l-1)(q^{l/2}-1),\quad g=l$ \par
\item[]  $l$ odd:  $p(q)=(q^2-1)^2(q^4-1)^2 \cdots (q^{l-1}-1)^2(q^{(l+1)/2}-1),\quad g=l$ \par
\end{itemize}

$D_l$:  $\check K=SO(l)\times SO(l)$,
 \begin{itemize}
\item[]   $l$ even: $p(q)=(q^2-1)^2(q^4-1)^2(q^6-1)^2...(q^{l-2}-1)^2(q^{l/2}-1)^2,\quad g=l$ \par
\item[]   $l$ odd: $p(q)=(q^2-1)^2(q^4-1)^2(q^6-1)^2...(q^{l-1}-1)^2,\quad g=l-1$ \par
\end{itemize}
    
$E_6$:  ${\rm Lie}(\check K)=\mathfrak{sp}(4)$,
\begin{itemize}
\item[]   $p(q)=(q^2-1)(q^4-1)(q^6-1)(q^8-1),\quad g=4$
\end{itemize}

 $E_7$:   ${\rm Lie}(\check K)={\mathfrak{su}}(8)$,
\begin{itemize}
\item[]  $p(q)=(q^2-1)(q^3-1)(q^4-1)(q^5-1)(q^6-1)(q^7-1)(q^8-1),\quad g=7$ 
\end{itemize}

$E_8$: ${\rm Lie}(\check K)={\mathfrak{so}}(16)$,
\begin{itemize}
\item[]  $p(q)=(q^2-1)(q^4-1)(q^6-1)(q^8-1)(q^{10}-1)(q^{12}-1)(q^{14}-1)(q^8-1),\quad g=8$ \par
\end{itemize}

$F_4$:  ${\rm Lie} (\check K)={\mathfrak{sp}}(1)\times {\mathfrak{sp}}(3)$,
 \begin{itemize}
 \item[]  $p(q)=(q^2-1)^2(q^4-1)(q^6-1),\quad g=4$   
\end{itemize}

$G_2$:  ${\rm Lie}(\check K)=\mathfrak{su}(2)\times {\mathfrak{su}}(2)$,
  \begin{itemize}
\item[]  $p(q)=(q^2-1)^2,\quad g=2$ 
\end{itemize}

 \end{Theorem}
In the factorization $p(q)= \Pi_{i=1}^g (q^{d_i}-1)$, 
 the $d_i$ are the degrees of the basic Weyl group invariant
polynomials for the compact Lie group $K$ (see \cite{carter:89}). The numbers $2d_i-1$ are the degrees
of the generators of the cohomology ring of $K$, i.e. $H^*(K;{\mathbb Q})\cong
\Lambda(x_1,\ldots,x_g)$ with deg$(x_i)=2d_i-1$. In this
sense the cohomology ring of $K$ over the rationals can be derived from
the structure of the blow-up points of the Toda lattice. This can be made
more explicit and also allows us to write cocycles representing
the generators of $K$ in terms of duals of the Bruhat cells.

\section{Zeros of Schur polynomials and $p(q)$}
\label{singularstructure}

We now show that
the singular structure of the Painlev\'e divisor ${\mathcal D}_0=\cup_{j=1}^l{\mathcal D}_j$ at
the point $p_o$ is related to the zeros of certain Schur polynomials. We have the following conjecture
for the nonperiodic case:

\begin{Conjecture}\label{tangentcone}
The degree $\eta(w_*)$ of the polynomial $p(q)$ is related to the multiplicity $d$ of the singularity of the
divisor given by ${\mathcal D}_0=\cup_{j=1}^l{\mathcal D}_j$ at the point $p_o$
where all the divisors ${\mathcal D}_j=\{\tau_j=0\}$ intersect. Namely the number $d$ is
given by the minimal degree of the product of the $\tau$-functions,
that is, the degree of the tangent cone given by (\ref{divF}).
Furthermore, the degree $\eta(w_*)$
also gives the number of {\it real} $t_1$ roots of the product of the Schur polynomials 
$S_{\lambda_k}(t_1,\ldots,t_l)$ for generic values of $t_2,\ldots, t_l$.  Here $\lambda_k,~k=1,\ldots,l$,
indicate the Young diagrams $\lambda_k$, and those Schur polynomials are associated to
the $\tau_k$-function of the nilpotent Toda lattices.
\end{Conjecture}
The first part of the conjecture can be verified directly for the Toda lattice associated with any Lie algebra of the classical type or type $G_2$
by counting the minimal degrees of Schur polynomials,
which are the leading terms of the $\tau$-functions near the point $p_o$.
The second part is verified for the cases with lower ranks (see below).

\subsection{Schur polynomials appearing in the nilpotent Toda lattices}
Let us first show that the $\tau$-functions in the semisimple case near the point $p_o$ can be approximated by
those in the nilpotent case (see also \cite{casian:04}). We here explain the case of $A_l$, and other cases of
$B,C,D$ and $G$ can be discussed in the similar way.
In the case of $A_l$, the $\tau_1$ function defined by $\tau_1=\langle g e_1,e_1\rangle$
with $g=\exp(\sum_{j=1}^l(L^0)^jt_j)$, i.e. (\ref{tau}) can be expressed by
\[\begin{array}{llll}
\tau_1(t_1,\ldots,t_l)&=&\displaystyle{\sum_{k=0}^{l}\rho_k\exp\left(\sum_{j=1}^l\lambda_kt_j\right)}\\
&=&\displaystyle{\sum_{n=0}^{\infty}h_n(t_1,\ldots,t_l)\left(\sum_{k=0}^l\lambda_k^n\rho_k\right)}\,,
\end{array}
\]
where $\lambda_k$ are the eigenvalues of the $(l+1)\times(l+1)$ Lax matrix $L^0$, and $\rho_k$ are constants
determined by $L^0$. The functions $h_k(t_1,\ldots,t_l)$ are the complete homogeneous symmetric functions given by
\begin{equation}
\label{schurp}
 h_k(t_1,\ldots,t_l)=\sum_{i_1+2i_2+\cdots+li_l=k}\frac{t_1^{i_1}t_2^{i_2}\cdots t_l^{i_l}}{
 i_1!i_2!\cdots i_l!}\,.
 \end{equation}
  We set $t=(0\ldots,0)$ to be the point $p_o$, that is,
\[
\tau_1(0,\ldots,0)=\tau_2(0,\ldots,0)=\cdots=\tau_l(0,\ldots,0)=0\,.
\]
This implies $(\partial^k\tau_1/\partial t_1^k)(0,\ldots,0)=0$ for $k=0,1,\ldots,l-1$, and $(\partial^l\tau_1/\partial t_1^l)(0,\ldots,0)=1$, from which we obtain
\[
\tau_1(t_1,\ldots,t_l)=\displaystyle{\sum_{n=l}^{\infty}h_n(t_1,\ldots,t_l)\left(\sum_{k=0}^l\lambda_k^n\rho_k\right)}\,,
\]
where $\rho_k$ are determined by
\[\begin{pmatrix}
1   &    1   &  \cdots  &  1  \\
\lambda_0&\lambda_1&\cdots&\lambda_l \\
\vdots & \vdots  &  \ddots  & \vdots  \\
\lambda_0^l &\lambda_1^l &\cdots &\lambda_l^l
\end{pmatrix}\begin{pmatrix} \rho_0 \\ \rho_1\\ \vdots \\ \rho_l \end{pmatrix}=
\begin{pmatrix} 0 \\ 0 \\ \vdots \\ 1 \end{pmatrix}
\]
Notice that in the nilpotent case (all $\lambda_k=0$), we have $\tau_1=h_{l}=S_{(l)}$
(recall that $\tau_1=\langle ge_{l+1},e_1\rangle$ with $g\in G^{C_0}$).
Now using the equation of the $\tau$-functions which is derived from (\ref{toda-lax}) and (\ref{ab}), i.e.
\begin{equation}\label{hirota}
\tau_k\frac{\partial^2\tau_k}{\partial t_1^2}-\left(\frac{\partial \tau_k}{\partial t_1}\right)^2=a_k^0\prod_{j\ne k}(\tau_j)^{-C_{k,j}}\,,
\end{equation}
we obtain 
\begin{equation}
\label{schur-tau}
\tau_k=(-1)^{\frac{k(k-1)}{2}}S_{(l-k+1,\ldots,l)}(t_1,\ldots,t_l)+ ({\rm higher~degree~terms})\,.
\end{equation}
where $S_{(i_1,\ldots,i_k)}$ with $1\le i_1<\cdots<i_k\le l$ is the Schur polynomial defined as
the Wronskian determinant with respect to the $t_1$-variable (see \cite{casian:04}),
\[
S_{(i_1,\ldots,i_k)}={\rm Wr}(h_{i_1},\ldots,h_{i_k})=\Big|\left(h_{i_{\alpha}-\beta+1}\right)_{1\le\alpha,\beta,\le k}\Big|\,.
\]
This implies that the $\tau$-functions in the semisimple case near the point $p_o$ can be
approximated by the Schur polynomials which are the $\tau$-functions
in the nilpotent case.
In particular, the multiplicity of the zero of $\tau_k$ is given by $k(l-k+1)$ which is the degree of
the Schur polynomial in (\ref{schur-tau}). Also notice that the $\tau$ functions are weighted homogeneous polynomials where the weight is defined by $k$ for $h_k$.

       The multiplicities of the $\tau$ functions for any Lie algebra ${\mathfrak g}$
can be obtained from (\ref{hirota}) (see also \cite{flaschka:91}): Let $\nu_k$ be the multiplicity of the zero for $\tau_k(t_1)$, that is, $\tau_k(t_1)\sim t_1^{\nu_k}$. Then from (\ref{hirota}), we have 
$2(\nu_k-1)=-\sum_{j\ne k}C_{k,j}m_j$ which leads to
\[
\nu_k=2\sum_{j=1}^l (C^{-1})_{k,j}\,,
\]
(see also Proposition 2.4 in \cite{flaschka:91}).
This number is also related to the total number of the root $\alpha_k$ in the positive root system $\Delta^+$ for the {\it Langlands dual} algebra $\check{\mathfrak g}$: To show this, we calculate $2\rho$,
the sum of all the positive roots in ${\mathfrak g}$,
 \[
 2\rho=\sum_{\alpha\in\Delta_+}\alpha=\sum_{k=1}^l n_k\alpha_k\,.
 \]
 Note that $\rho$ is also given by the sum of the fundamental weights, $\rho=\sum_{k=1}^l\omega_k$.
  Then using the relation $\alpha_i=\sum_{j=1}^lC_{i,j}\omega_j$, the number $n_k$ is given by
 \[
 n_k=2\sum_{j=1}^l(C^{-1})_{j,k}\,.
 \]
 Thus the multiplicity $\nu_k$ for the $\tau_k$ function associated to ${\mathfrak g}$ is related to the number
 $n_k$ for the Langlands dual algebra $\check {\mathfrak g}$ whose Cartan matrix is given by
 the transpose of that for ${\mathfrak g}$. This duality leads to the appearance of
 the Langlands dual in Theorems \ref{mainT} and \ref{Kpq} (see \cite{casian:05} for the details). The multiplicity of the product of the $\tau$ functions at the point $p_o$ is then given by (see also \cite{adler:91,flaschka:91})
 \[
 |2\rho|=\sum_{k=1}^l\nu_k=\sum_{k=1}^ln_k=2\sum_{i,j=1}^l\left(C^{-1}\right)_{i,j}\,.
 \]

In the general case, the $\tau$-functions for the nilpotent
Toda lattices are obtained from (\ref{tau}). Then one finds explicit forms of the highest weight vectors and
the companion matrix (the regular nilpotent element) for each algebra. The following is the result
based on the $\tau$-functions of the nilpotent Toda lattices (see \cite{casian:04} for the nilpotent Toda lattices in general):

\bigskip
\noindent
{\bf $A_l$-Toda lattices:}
The $\tau$-functions
are given by the Schur polynomials in (\ref{schur-tau}), i.e.
\[
\tau_{k}(t_1,\ldots,t_l)=(-1)^{\frac{k(k-1)}{2}}S_{(l-k+1,\ldots,l)}(t_1,\ldots,t_l)\quad k=1,\ldots,l\,.
\]
The minimal degrees of those Schur polynomials are given as follows:
\begin{itemize}
\item{} For $l$ even, the minimal degrees of $\tau_j$, $j=1,\ldots,l$, are given by
\[
1,~2,~\ldots~,~\frac{l}{2},~\frac{l}{2},~,\ldots~,~2,~1,
\]
(e.g. $\tau_1\sim t_l, ~\tau_2\sim t_{l-1}^2$ and $\tau_l\sim t_l$).
 The sum of those degrees then gives $d=\frac{l(l+2)}{4}$.
\item{} For $l$ odd, the minimal degrees are
\[
1,~2,~\ldots~,~\frac{l-1}{2},~\frac{l+1}{2},~\frac{l-1}{2},~\ldots~,~2,~1.
\]
 from which we have $d=\left(\frac{l+1}{2}\right)^2$.
\end{itemize}
For example, in the case of $l=2$, we have $\tau_1=S_{(2)}=t_2+\frac{t_1^2}{2},\,\tau_2=-S_{(1,2)}=t_2-\frac{t_1^2}{2}$. Thus we have the degree four polynomial for $F(t_1,t_2)=\tau_1\tau_2(t_1,t_2)$,
and the minimal degree is two which is the number of blow-ups $\eta(w_*)$. Note here
that $|2\rho|=4$ is the number of complex roots.  The second part of Conjecture \ref{tangentcone} then states  that
$F(t_1,t_2)=-S_{(1)}S_{(1,2)}(t_1,t_2)$ has two {\it real} roots in $t_1$ for a generic value of $t_2=$constant (i.e. $t_2\ne 0$). One should also note that the sum of the minimal degrees 
of the pair $\tau_k$ and $\tau_{l-k+1}$ is equal to the degree $d_k$ in Theorem \ref{Kpq}.

\bigskip
\noindent
{\bf $B_l$-Toda lattice:} The $\tau$-functions are given by
\[
\tau_k(t_1,t_3,\ldots,t_{2l-1})={\rm Wr}(h_{2l},\ldots,h_{2l-k+1})\, \quad k=1,\ldots,l-1
\]
and 
\[\tau_l(t_1,t_3,\ldots,t_{2l-1})=\sqrt{|{\rm Wr}(h_{2l},\ldots,h_{l+1})|}\,,
\]
where $h_k$ are given in (\ref{schurp}) with $t_{2k}=0$ for all even parameters (see \cite{casian:04}).
Note the indices of the flow parameters, $2k-1$, $k=1,\ldots,l$, are given by the exponents $m_k$ of the root system $\Delta$.
Note that the $\tau_l$ is given by the pfaffian related to the Schur $Q$-function. More precisely, we
have $\tau_l(t_1,t_3,\ldots,t_{2l-1})\sim S_{(1,3,\ldots,2l-1)}(t_1/2,t_3/2,\ldots,t_{2l-1}/2)$ which is 
the Schur $Q$-polynomial, $Q_{(l,l-1,\ldots,1)}(t_1,t_3,\ldots,t_{2l-1})$ \cite{fulton:98}.
Then we have:
\begin{itemize}
\item{} For $l$ even, the minimal degrees are given by
\[
2,~2,~4,~4,~\ldots~,~l-2,~l-2,~l,~\frac{l}{2}\,.
\]
(e.g. $\tau_1\sim t_1t_{2l-1},~\tau_2\sim t_{2l-1}^2$ and $\tau_l\sim (t_{l+1})^{l/2}$).
The degree $d$ is then given by $d=\frac{l(l+1)}{2}$.
\item{} For $l$ odd, the minimal degrees are
\[
2,~2,~4,~4,~\ldots~,~l-1,~l-1,~\frac{l+1}{2}\,,
\]
 from which we have $d=\frac{l(l+1)}{2}$.
\end{itemize}
For the case $B_2$, we have $\tau_1=S_{(4)}=t_1(t_3+\frac{t_1^3}{24})$ and $\tau_2=\sqrt{S_{(3,4)}}=
t_3-\frac{t_1^3}{12}$. The degree of $F=\tau_1\tau_2$ in $t_1$ is seven which is the hight $|2\rho|$, and
the minimal degree is three which is the number of blow-ups $\eta(w_*)$ and also the number of real roots of $F(t_1,t_3)$ in $t_1$
for $t_3\ne 0$. One should note that the minimal degree of $\tau_k$ also gives the degree $d_k$
in Theorem \ref{Kpq} for the Langlands dual algebra $C_l$.

\bigskip
\noindent
{\bf $C_l$-Toda lattice:} The $\tau$-functions are
\[
\tau_k(t_1,t_3,\ldots,t_{2l-1})={\rm Wr}(h_{2l-1},\ldots,h_{2l-k})\,\quad k=1,\ldots,l\,.
\]
Again one takes $t_{2k}=0$ for $h_n$. Then the minimal degrees are given by
\[
1,~2,~3,~\ldots~,~l-1,~l\,.
\]
This gives $d=\frac{l(l+1)}{2}$ (which is the same as $B_l$-case).
For the case $C_2$, we have $\tau_1=S_{(3)}=t_3+\frac{t_1^3}{6}$ and $\tau_2=-S_{(2,3)}=t_1(t_3-\frac{t_1^3}{12})$.
The product $F=\tau_1\tau_2$ has the degree seven which is the height $|2\rho|$, and
the minimal degree is three which is $\eta(w_*)$ and also the number of real roots of $F(t_1,t_3)$ in $t_1$
for $t_3\ne 0$.  The minimal degree of $\tau_k$ then gives $d_k$ for the Langlands dual
algebra $B_l$ in Theorem \ref{Kpq}.

\bigskip
\noindent
{\bf $D_l$-Toda lattice:} The $\tau$-functions are given as follows:
\begin{itemize}
\item{} For $l$ even, they are given by, for $k=1,\ldots,l-2$,
\[
\tau_k(t_1,t_3,\ldots,t_{2l-3},s)={\rm Wr}(sh_{l-1}+2h_{2l-2},sh_{l-2}+2h_{2l-3},\ldots,sh_{l-k}+2h_{2l-1-k})\,.
\]
The $\tau_{l-1}$ and $\tau_l$ are given by
\[
[\tau_{l-1}\cdot\tau_l](t_1,t_3,\ldots,t_{2l-3},s)={\rm Wr}(sh_{l-1}+2h_{2l-2},\ldots,sh_{2}+2h_l)\,,
\]
and
\[
(\tau_l(t_1,t_3,\ldots,t_{2l-3},s))^2=\pm\left|\begin{matrix}
sh_{l-1}+2h_{2l-2} & sh_{l-2}+2h_{2l-3}&\cdots &sh_{1}+2h_{l}&s+h_{l-1}\\
sh_{l-2}+2h_{2l-3} &sh_{l-3}+2h_{2l-4}&\cdots & s+2h_{l-1}  &  h_{l-2} \\
\vdots&\vdots  & \ddots  &  \vdots &\vdots  \\
sh_1+2h_l & s+2h_{l-1} &\cdots & 2h_2 & h_1\\
s+h_{l-1} & h_{l-2} &\cdots & h_1 & 0
\end{matrix}\right|
\]
Here the even parameters are all zero $t_{2k}=0$, and $s$ is a flow parameter
associated with the Chevalley invariant with the degree $l$, i.e. the exponent $m_l=l-1$.
The exponents $m_k$ of the root system of $D$-type are
given by $(1, 3,\ldots,2l-3, l-1)$ (see e.g. \cite{carter:89}).
Counting the minimal degrees of the $\tau$-functions, we have
\[
2,~2,~4,~4,~\ldots~,~l-2,~l-2,~\frac{l}{2},~\frac{l}{2}\,,
\]
(e.g. $\tau_1\sim st_{l-1},~\tau_2\sim t_{2l-3}^2$ and $\tau_l\sim s^{l/2}$). Then we have $d=\frac{l^2}{2}$. The minimal degree of $\tau_k$ then gives the degree $d_k$ in
Theorem \ref{Kpq}.
\item{} For $l$ odd, the $\tau$-functions are given by, for $k=1,\ldots,l-2$,
\[
\tau_k(t_1,t_3,\ldots,t_{2l-3},s)={\rm Wr}(s^2+2h_{2l-2},2h_{2l-3},\ldots,2h_{2l-1-k})\,.
\]
The last two $\tau$-functions are 
\[
[\tau_{l-1}\cdot\tau_l](t_1,t_3,\ldots,t_{2l-3},s)={\rm Wr}(s^2+2h_{2l-2},2h_{2l-3},\ldots,2h_{l})\,,
\]
and 
\[
(\tau_l(t_1,t_3,\ldots,t_{2l-3},s))^2=\pm\left|\begin{matrix}
s^2+2h_{2l-2} & 2h_{2l-3}&\cdots &2h_{l}&s+h_{l-1}\\
2h_{2l-3} &2h_{2l-4}&\cdots & 2h_{l-1}  &  h_{l-2} \\
\vdots&\vdots  & \ddots  &  \vdots &\vdots  \\
2h_l & 2h_{l-1} &\cdots & 2h_2 & h_1\\
s+h_{l-1} & h_{l-2} &\cdots & h_1 & 1
\end{matrix}\right|
\]
Now the minimal degrees of the $\tau$-functions are
\[
2,~2,~\ldots~,~l-3,~l-3,~l-1,~\frac{l-1}{2},~\frac{l-1}{2}\,.
\]
Then we have $d=\frac{l^2-1}{2}$.
\end{itemize}

\bigskip
\noindent
{\bf $G_2$-Toda lattice:} We have
\[
\tau_1=S_{(6)}=t_1\left(t_5+\frac{t_1^5}{720}\right),\quad \tau_2=-S_{(5,6)}=t_5^2-\frac{t_5t_1^5}{40}+\frac{t_1^{10}}{86400}\,.
\]
Here the flow parameters are only $t_1$ and $t_5$ and all others take zero, i.e. $t_2=t_3=t_4=0$.
Those indices $1,5$ in the nonzero parameters are the exponents of the Weyl group of $G$-type.
The degree of $F=\tau_1\tau_2$ is $16$ $(=|2\rho|)$, and the minimal degree
is four which is $\eta(w_*)$, and this also gives the number of real roots of $F(t_1,t_5)$ in $t_1$
for $t_5\ne 0$. Each minimal degree of $\tau_k$ is also $d_k$ in Theorem \ref{Kpq}.

\subsection{The degree of the tangent cone at $p_o$ and $p(q)$}

We can now state the following Proposition from those computations:

\begin{Proposition} \label{direct} Let  ${\mathfrak g}$ be a split semisimple Lie algebra not containing factors of type $E$ or $F$. The degree $d$ of the tangent cone at $p_o$ of the Painlev\'e divisor
${\mathcal D}_0=\cup_{j=1}^l{\mathcal D}_j$, which is given as the minimal degrees
of Schur polynomials  in $(\ref{schur-tau})$ (see (\ref{divF})),  is the dimension of any Borel subalgebra of
${\rm Lie}( \check K  ({\mathbb C}))$. Moreover, $d$ is also the degree of the polynomial
$\tilde p(q)=q^{-r}|\check K({\mathbb F}_q)|$. 
\end{Proposition}

 Then from Theorem \ref{Kpq},  we obtain the following Proposition (the proof
 can be found in \cite{casian:05}):

 \begin{Proposition} 
 \label{totalBnumber}
 The number $\eta(w_*)={\rm deg}(p(q))$ satisfies
 $\eta(w_*)=d_1+\cdots +d_l$. Moreover we have $\eta(w_*)=d$ for any semisimple Lie algebra not
 containing factors of type $E$ or $F$.
 We have the following formulas for  $\eta(w_*)$. The number $\eta(w_*)$ is, in each case,  the complex dimension of any Borel subalgebra of  ${\rm Lie}(\check K({\mathbb C}))$, 
 \begin{itemize}
\item[]  $A_l$:  $\eta(w_*) =\frac{l(l+2)}{4}$ if $l$ is even;$\quad\eta(w_*)=\frac{(l+1)^2}{4}$ if $l$ is odd,
\item[]  $B_l$  or $C_l$:  $\eta(w_*)=\frac{l(l+1)}{2}$,
\item[]  $D_l$:  $\eta(w_*)=\frac{l^2}{2}$ if $l$ is even;$\quad\eta(w_*)=\frac{l^2-1}{2}$ if $l$ is odd,
\item[]  $E_l$:  $\eta(w_*)=20$ if $l=6$;$\quad\eta(w_*)=35$ if $l=7$;$\quad\eta(w_*)=64$ if $l=8$,
\item[]  $F_4$:  $\eta(w_*)=14$,
\item[]  $G_2$:  $\eta(w_*)=4$.

\end{itemize}
 \end{Proposition}

 We especially note that $d$ gives the multiplicity of a singularity at $p_o$, and  $\eta(w_*)$ is defined differently as the maximal number of blow-ups encountered along the Toda flow, counted along one dimensional subsystems. The polynomials $p(q)$, the alternating sum of the blow-ups, are shown to agree with  the polynomials  $\tilde p(q)$. As was noticed, the minimal degree
for each $\tau$-functions is related to each degree $d_i$ of the basic $W$-invariant
polynomial of the Chevalley group $\check K$. For example, the
degree $d_i$ and the minimal degree of $\tau_i$-function are the same for the cases having
the same ranks, $l={\rm rank}(\mathfrak g)={\rm rank}(K)$, i.e.
the cases of $B,\,C$ and $D_l$ with $l$ even. We did not compute  the exact relation between those
degrees, but we expect that each degree $d_i$ is related to the number of real intersection points
on the tangent cone $V=\{F_d=0\}$ defined in (\ref{divF}) with a linear line
corresponding to the $t_1$-flow of the Toda lattice. This may be stated as
\[
d={\rm deg}(p(q))=\underset{c\in U}{\rm Max}\left|\{{\mathcal D}_0\cap L_c~|~
{\rm transversal ~intersection}\,\}\right|
\]
where $U\subset {\mathbb R}^{l-1}$ is a neighborhood of t=0, and $L_c:=\{(t_1,c_2,\ldots,c_l)\,:\,c=(c_2,\ldots,c_l)\in U\}$
is the linear line of $t_1$-flow. This statement is equivalent to the second part of Conjecture \ref{tangentcone}, that is, the number of real $t_1$-roots of the product $F=\prod_{j=1}^{l}\tau_j$
in the nilpotent limit.

\bibliographystyle{amsalpha}

\end{document}